\documentclass[10pt,a4paper]{article}
\usepackage{mathrsfs}
\usepackage{amsfonts}
\usepackage{}

\usepackage[leqno]{amsmath}
\usepackage{graphicx}
\usepackage{latexsym}
\usepackage{amsmath,amsfonts,amssymb,amsthm,mathrsfs,euscript,makeidx,color}
\usepackage{enumerate}

\usepackage[colorlinks,linkcolor=blue,anchorcolor=green,citecolor=red]{hyperref}

\def\sqr#1#2{{\vcenter{\vbox{\hrule height.#2pt
              \hbox{\vrule width.#2pt height#1pt \kern#1pt \vrule width.#2pt}
              \hrule height.#2pt}}}}
\def\signed #1{{\unskip\nobreak\hfil\penalty50
              \hskip2em\hbox{}\nobreak\hfil#1
              \parfillskip=0pt \finalhyphendemerits=0 \par}}
\def\endpf{\signed {$\sqr69$}}

\def\5n{\negthinspace \negthinspace \negthinspace \negthinspace \negthinspace }
\def\4n{\negthinspace \negthinspace \negthinspace \negthinspace }
\def\3n{\negthinspace \negthinspace \negthinspace }
\def\2n{\negthinspace \negthinspace }
\def\1n{\negthinspace }

\def\dbE{\mathbb{E}}
\def\dbF{\mathbb{F}}

\def\dbM{\mathbb{M}}

\def\dbP{\mathbb{P}}

\def\dbR{\mathbb{R}}

\def\dbX{\mathbb{X}}
\def\dbY{\mathbb{Y}}

\def\sD{\mathscr{D}}

\def\sL{\mathscr{L}}

\def\sR{\mathscr{R}}
\def\sS{\mathscr{S}}
\def\sT{\mathscr{T}}

\def\cA{{\cal A}}
\def\cB{{\cal B}}
\def\cC{{\cal C}}
\def\cD{{\cal D}}

\def\cF{{\cal F}}
\def\cG{{\cal G}}

\def\cL{{\cal L}}
\def\cM{{\cal M}}

\def\cQ{{\cal Q}}
\def\cR{{\cal R}}
\def\cS{{\cal S}}

\def\cU{{\cal U}}

\def\cX{{\cal X}}
\def\cY{{\cal Y}}
\def\cZ{{\cal Z}}

\def\BA{{\bf A}}
\def\BB{{\bf B}}
\def\BC{{\bf C}}
\def\BD{{\bf D}}

\def\BG{{\bf G}}

\def\BQ{{\bf Q}}
\def\BR{{\bf R}}
\def\BS{{\bf S}}

\def\BX{{\bf X}}

\def\Bg{{\bf g}}

\def\Bq{{\bf q}}

\def\Bu{{\bf u}}

\def\Bxi{{\xi\3n\xi\3n\xi\3n\xi\3n\xi}}
\def\Brho{{\rho\3n\rho\3n\rho\3n\rho\3n\rho}}

\def\BBeta{{\eta\3n\eta\3n\eta}}

\def\ae{\hbox{\rm a.e.}}
\def\as{\hbox{\rm a.s.}}

\def\ds{\displaystyle}

\def\ns{\noalign{\ss}}
\def\no{\noindent}

\def\ss{\smallskip}
\def\ms{\medskip}
\def\bs{\bigskip}
\def\q{\quad}
\def\qq{\qquad}
\def\hb{\hbox}

\def\({\Big (}
\def\){\Big )}
\def\[{\Big[}
\def\]{\Big]}

\def\lan{\langle}
\def\ran{\rangle}

\def\blan{\big\langle}
\def\bran{\big\rangle}

\def\rf{\eqref}

\def\a{\alpha}

\def\g{\gamma}
\def\d{\delta}
\def\e{\varepsilon}

\def\l{\lambda}
\def\m{\mu}
\def\n{\nu}
\def\si{\sigma}
\def\t{\tau}
\def\f{\varphi}
\def\th{\theta}
\def\o{\omega}

\def\i{\infty}
\def\G{\Gamma}

\def\F{\Phi}
\def\O{\Omega}

\def\esssup{\mathop{\rm esssup}}

\def\h{\widehat}
\def\wt{\widetilde}

\def\cd{\cdot}
\def\cds{\cdots}

\def\les{\leqslant}
\def\ges{\geqslant}
\def\bde{\begin{definition}\label}
\def\ede{\end{definition}}
\def\be{\begin{equation}}
\def\bel{\begin{equation}\label}
\def\ee{\end{equation}}
\def\bt{\begin{theorem}\label}
\def\et{\end{theorem}}
\def\bc{\begin{corollary}\label}
\def\ec{\end{corollary}}
\def\bl{\begin{lemma}\label}
\def\el{\end{lemma}}
\def\bp{\begin{proposition}\label}
\def\ep{\end{proposition}}
\def\bas{\begin{assumption}\label}
\def\eas{\end{assumption}}
\def\br{\begin{remark}\label}
\def\er{\end{remark}}
\def\bex{\begin{example}\label}
\def\ex{\end{example}}
\def\ba{\begin{array}}
\def\ea{\end{array}}
\def\ben{\begin{enumerate}}
\def\een{\end{enumerate}}

\def\square#1{\vbox{\hrule\hbox{\vrule height#1%
     \kern#1\vrule}\hrule}}
\def\rectangle#1#2{\vbox{\hrule\hbox{\vrule height#1%
     \kern#2\vrule}\hrule}}

\font\tenbb=msbm10 \font\sevenbb=msbm7 \font\fivebb=msbm5

\newfam\bbfam
\scriptscriptfont\bbfam=\fivebb \textfont\bbfam=\tenbb
\scriptfont\bbfam=\sevenbb

\newtheorem{theorem}{\indent Theorem}[section]
\newtheorem{definition}[theorem]{\indent Definition}
\newtheorem{proposition}[theorem]{\indent Proposition}
\newtheorem{corollary}[theorem]{\indent Corollary}
\newtheorem{lemma}[theorem]{\indent Lemma}
\newtheorem{remark}[theorem]{\indent Remark}
\newtheorem{example}[theorem]{\indent Example}

\newtheorem{assumption}[theorem]{\indent Assumption}

\makeatletter
   
   \@addtoreset{equation}{section}
\makeatother

\begin{document}

\title{\bf Linear Quadratic Stochastic Optimal Control\\ Problems with Operator Coefficients:\\ Open-Loop Solutions\thanks{This work is supported in part
by NSF Grant DMS-1406776, the National Natural Science Foundation of China (11471192,
11401091,11571203), the Nature Science Foundation of Shandong Province (JQ201401), the Fundamental
Research Funds of Shandong University (2017JC016), and the Fundamental Research
Funds for the Central Universities (2412017FZ008).}}

\author{Qingmeng Wei\footnote{School of Mathematics and Statistics, Northeast Normal University, Changchun 130024, China; email: {\tt weiqm100@nenu.} {\tt edu.cn}},~~
Jiongmin Yong\footnote{Department of Mathematics, University of
Central Florida, Orlando, FL 32816, USA; email: {\tt
jiongmin.yong@ucf.edu}},~~and~~Zhiyong Yu\footnote{Corresponding
author, School of Mathematics, Shandong University, Jinan 250100,
China; email: {\tt yuzhiyong@sdu.edu.cn}}}

\maketitle

\no\bf Abstract: \rm An optimal control problem is considered for linear stochastic differential equations with quadratic cost functional. The coefficients of the state equation and the weights
in the cost functional are bounded operators on the spaces of square integrable random variables. The main motivation of our study is linear quadratic (LQ, for short) optimal control problems for
mean-field stochastic differential equations. Open-loop solvability of the problem is characterized as the solvability of a system of linear coupled forward-backward stochastic differential equations
(FBSDE, for short) with operator coefficients, together with a convexity condition for the cost functional. Under proper conditions, the well-posedness of such an FBSDE, which leads to the existence of an open-loop optimal control, is established. Finally, as applications of our main results, a general mean-field LQ control problem and a concrete mean-variance portfolio selection problem in
the open-loop case are solved.

\ms

\no\bf Keywords: \rm linear stochastic differential equation with
operator coefficients, open-loop solvability, forward-backward
stochastic differential equations, mean-field linear quadratic
control problem, mean-variance portfolio selection.

\ms

\no\bf AMS Mathematics Subject Classification. \rm 93E20, 91A23, 49N70, 49N10.

\section{Introduction}

Since the pioneer works of Lasry--Lions \cite{Lasry-Lions 2006-1,Lasry-Lions 2006-2,Lasry-Lions 2007}, the mean-field type stochastic differential equations, stochastic optimal controls,
stochastic differential games, and their applications received extensive researchers' attention in recent years (\cite{Buckdahn-Djehiche-Li-Peng 2009,Buckdahn-Li-Peng 2009,Bensoussan-Frehse-Yam 2013,Carmona-Delarue 2013,Yong 2013,Bensoussan-Yam-Zhang 2015,Carmona-Delarue
2015,Bensoussan-Sung-Yam-Yung 2016}). This paper is mainly motivated by the mean-field linear quadratic (LQ, for short) stochastic optimal control problems, and this motivation will be explained in
detail in the next section.

\ms

When we consider mean-field LQ control problems with random coefficients, the terms of the following forms
\bel{AX}A(s)X(s),\q\bar A(s)\dbE\big[\wt A(s)X(s)\big],\q\hb{etc.}\ee
could appear in the state equations, where $s$ is the time variable, $A(\cd)$, $\bar A(\cd)$, $\wt A(\cd)$ are suitable matrix-valued stochastic processes, and $X(\cd)$ denotes the state
process. On the other hand, the terms of the forms
\bel{Sec1_Items_MF_Cost}\lan Q(s)X(s),X(s)\ran,\q\big\lan\bar Q(s)X(s),\dbE\big[\wt Q(s)X(s)\big] \big\ran,\q\big\lan\wt Q(s)\dbE\big[\wt Q(s)X(s)\big],\dbE\big[\wt Q(s)X(s)\big]\big\ran,\q\hb{etc.}\ee
could appear in the cost functionals. It is a bit complicated to treat them case by case. This inspires us to introduce a universal framework to deal with them uniformly.

\ms

In order to further reveal the advantages of the framework of operators, let us  see a simple example. It is well known that mean-variance portfolio selection problems are very important in mathematical finance. In a natural viewpoint, mean-variance problems can be regarded as special examples of mean-field LQ problems (see \cite{Andersson-Djehiche 2011,Bjork-Murgoci-Zhou 2014, Pham-Wei 2017} or Subsection 4.1 in the present paper). In mean-variance problems, the variance of some
random variable is included in the cost functionals. For example, the variance of the terminal state $X(T)$ reads:
\bel{}\hb{Var}\big( X(T)\big)=\dbE\[X(T)^2-2X(T)\dbE\big[X(T)\big]+\big(\dbE\big[X(T)\big]\big)^2\]=\dbE\[X(T)^2
-\big(\dbE\big[X(T)\big]\big)^2\].\ee
Therefore, in the cost functional, if we allow the variances of the state and/or control to appear, then the expectations will appear quadratically. Having this in mind, it will be very natural to include terms like those in \eqref{Sec1_Items_MF_Cost} in the cost functional. We will present some concrete examples in a later section.

\ms

Now we introduce our framework. Let $(\O,\cF,\dbF,\dbP)$ be a complete filtered probability space on which a standard one-dimensional Brownian motion $\{W(t),t\ges0\}$ is defined such that $\dbF=\{\cF_t\}_{t\ges0}$ is the natural filtration of $W(\cd)$ augmented by all the $\dbP$-null sets in $\cF$. Consider the following controlled linear (forward) stochastic differential equation (FSDE, for short) on $[t,T]$:
\bel{state}\left\{\2n\ba{ll}
\ds dX(s)=\big[\cA(s)X(s)+\cB(s)u(s)+b(s)\big]ds+\big[\cC(s)X(s)+\cD(s)u(s)+\si(s)\big]dW(s),\q s\in[t,T],\\
\ns\ds X(t)= x.\ea\right.\ee
In the above, $X(\cd)$ is called the {\it state process} taking
values in the $n$-dimensional Euclidean space $\dbR^n$; $u(\cd)$ is
called the {\it control process} taking values in $\dbR^m$; $(t,x)$
is called an {\it initial pair} with $t\in[0,T)$ and $x$ being a square
integrable $\dbR^n$-valued $\cF_t$-measurable random variable; $b(\cd)$ and
$\si(\cd)$ are called {\it non-homogeneous terms}. To explain the
coefficients of the system, we first recall the following spaces:
For any $t\in[0,T]$,
$$\ba{ll}
\ns\ds L^2_{\cF_t}(\O;\dbR^n)=\big\{\xi:\O\to \mathbb
R^n\bigm|\xi\hb{ is $\cF_t$-measurable and $\|\xi\|_2\equiv\big(\dbE|\xi|^2\big)^{1\over2}<\i$
}\big\},\q L^2(\O;\dbR^n)=L^2_{\cF_T}(\O;\dbR^n),\\ [2mm]
\ns\ds L^2_\dbF(t,T;\dbR^n)=\Big\{\f:[t,T]\times\O\to \mathbb
R^n\bigm|\f(\cd)\hb{ is $\dbF$-progressively measurable, }
\dbE\int_t^T|\f(s)|^2ds<\i\Big\},\\
%
%
\ns\ds L^2_{\dbF}(\O;L^1(t,T;\dbR^n))=\Big\{\f:[t,T]\times\O\to\dbR^n\bigm|\f(\cd)\hb{ is $\dbF$-progressively measurable, } \Big\| \int_t^T |\f(r)|dr\Big\|_2 <\infty \Big\},\\
\ns\ds L^2_\dbF(\O;C([t,T];\dbR^n))=\Big\{\f:[t,T]\times\O\to\dbR^n\bigm|\f(\cd)\hb{ is $\dbF$-progressively measurable}, s\mapsto\f(s,\o)\hb{ is continuous}\\
\ns\ds\qq\qq\qq\qq\qq\qq\qq\qq\qq\qq\qq\qq\qq  \hb{almost surely, }
\mathbb E\[\sup_{s\in[t,T]}|\f(s)|^2\]<\i\Big\}.\ea$$
%
For any Banach spaces $\dbX$ and $\dbY$, we let $\sL(\dbX;\dbY)$ be the set of all linear bounded operators from $\dbX$ to $\dbY$, and denote $\sL(\dbX;\dbX)=\sL(\dbX)$. Also, when $\dbX$ is a Hilbert
space, we let $\sS(\dbX)$ be the set of all bounded self-adjoint operators on $\dbX$. In the state equation \rf{state}, we assume that
\bel{cA}\cA(s),\cC(s)\in\sL\(L^2_{\cF_s}(\O;\dbR^n)\),\q\cB(s),\cD(s)\in\sL\(L^2_{\cF_s}(\O;\dbR^m);
L^2_{\cF_s}(\O;\dbR^n)\),\qq\forall s\in[0,T].\ee
More precisely, for example, for any $\xi\in L^2_{\cF_s}(\O;\dbR^n)$, $\cA(s)\xi\in L^2_{\cF_s}(\O;\dbR^n)$, with
$$\(\dbE|\cA(s)\xi|^2\)^{1\over2}\equiv\|\cA(s)\xi\|_2\les\|\cA(s)\|\,\|\xi\|_2,$$
where
$$\|\cA(s)\|=\sup\Big\{\|\cA(s)\xi\|_2\bigm|\xi\in L^2_{\cF_s}(\O;\dbR^n),~\|\xi\|_2=1\Big\}.$$
Some additional conditions will be assumed for $\|\cA(\cd)\|$, $\|\cC(\cd)\|$, $\|\cB(\cd)\|$ and $\|\cD(\cd)\|$ later. In what follows, the set of all initial pairs is denoted by
$$\sD=\Big\{(t,x)\bigm|t\in[0,T],~x\in L^2_{\cF_t}(\O;\dbR^n)\Big\},$$
and the set of all {\it admissible controls} on $[t,T]$ is denoted by $\cU[t,T]=L^2_\dbF(t,T;\dbR^m).$

\ms

One can show that under certain conditions, for any initial pair $(t,x)\in\sD$ and control $u(\cd)\in \cU[t,T]$, the state equation \rf{state} admits a unique strong solution $X(\cd)\equiv X(\cd\,;t,x,u(\cd))\in L^2_\dbF(\O;C([t,T];\dbR^n))$. The performance of the control process is measured by the following cost functional:
\bel{cost}\ba{ll}
\ns\ds J(t,x;u(\cd))=\dbE\[\lan\cG X(T),X(T)\ran+2\lan g,X(T)\ran+\int_t^T\(\lan\cQ(s)X(s),X(s)\ran+2\lan\cS(s)X(s),u(s)\ran\\
\ns\ds\qq\qq\qq\qq\qq\qq\qq\qq\qq+\lan\cR(s)u(s),u(s)\ran+2\lan q(s),X(s)\ran+2\lan\rho(s),u(s)\ran\)ds\],\ea\ee
where
\bel{cG}\ba{ll}
\ns\ds\cG\in\sS\(L^2(\O;\dbR^n)\),\qq\cQ(s)\in\sS\(L^2_{\cF_s}(\O;\dbR^n)\),\\
\ns\ds\cS(s)\in\sL\(L^2_{\cF_s}(\O;\dbR^n);L^2_{\cF_s}(\O;\dbR^m)\),\qq\cR(s)\in\sS\(L^2_{\cF_s}
(\O;\dbR^m)\),\qq\forall s\in[0,T],\ea\ee
with certain additional conditions, and $g\in L^2(\O;\dbR^n)$,
$q(\cd)\in L^2_\dbF(\O;L^1(0,T;\dbR^n))$, $\rho(\cd)\in
L^2_\dbF(0,T;\dbR^m)$.

\ms

Our optimal control problem can be stated as follows.

\ms

\bf Problem (OLQ). \rm For given $(t,x)\in\sD$, find a $\bar u(\cd)\in\cU[t,T]$, called an {\it open-loop optimal control}, such that
\bel{inf J}J(t,x;\bar u(\cd))=\inf_{u(\cd)\in\cU[t,T]}J(t,x;u(\cd)).\ee

\ms

The above Problem (OLQ) clearly includes the classical stochastic linear quadratic (LQ, for short) optimal control problem for which all the coefficients and quadratic weighting operators in the cost
functional are matrix-valued processes (\cite{Tang 2003, Tang 2015}). On the other hand, by allowing the coefficients of the state equation and the quadratic weighting operators in the cost functional to be linear bounded operators between Hilbert spaces of square integrable random variables, our problem will cover stochastic LQ optimal control problem for mean-field FSDEs with cost functionals also involving mean-field terms (such problems are referred to as MF-LQ problems). In \cite{Yong 2013}, for a simple MF-LQ problem (with deterministic coefficients), under proper conditions, optimal control is obtained via the solution to a system of Riccati equations. See \cite{Huang-Li-Yong 2015,Li-Sun-Yong 2016,Yong 2016} for some follow-up works.

\ms

There are a plenty of literatures on the classical LQ optimal control problems and two-person differential games. See \cite{Schmitendorf 1970, Lukes-Russell 1971, Eisele 1982, Chen-Li-Zhou 1998, Chen-Yong 2001, Ait Rami-Moore-Zhou 2001,  Tang 2003, Zhang 2005, McAsey-Mou 2006, Sun-Yong 2014, Sun-Yong 2016} for deterministic coefficient cases, and \cite{Yu 2013, Huang-Yu 2014, Du 2015, Tang
2015} for random coefficient cases. In \cite{Sun-Yong 2014} (see also, \cite{Sun-Yong 2016}) open-loop and closed-loop solvabilities/saddle points were introduced, and the following interesting equivalent relations were established for LQ optimal control problems with deterministic coefficients: The open-loop solvability of the LQ problem is equivalent to the solvability of an FBSDE, and the closed-loop solvability of the LQ problem is equivalent to the solvability of the corresponding Riccati equation. For two-person differential games, similar results are also valid. In the current paper, we focus on the open-loop solvability of our Problem (OLQ). The studies of closed-loop case and differential game problems will be carried out in our future publications. For the solvability of FBSDEs or Riccati equations arising in the classical LQ stochastic optimal control problems and stochastic differential game problems, one is referred to \cite{Hamadene 1998, Peng-Wu 1999, Yu 2015}. We may regard the current work as a continuation of \cite{Yong 2013, Huang-Li-Yong 2015, Li-Sun-Yong 2016} and \cite{Sun-Yong 2014, Sun-Yong 2016}.

\ms

Several contributions have been made in this work. Firstly, all the involved coefficients are operator-valued processes or operator-valued variables. The appearance of the operator coefficients in the state equation and the cost functional prompts us to develop some new  methods and techniques. Actually, our results on the FSDEs and backward stochastic differential equations (BSDEs, for short) with operator coefficients are of independent interests themselves. Secondly, under a convexity condition, we establish the equivalence between the open-loop solvability of Problem (OLQ) and the well-posedness of a coupled FBSDE with operator coefficients. Thirdly, under some conditions, the well-posedness of the relevant FBSDE with operator coefficients is established by the method of continuation (which is original introduced in \cite{Hu-Peng 1995}). Fourthly, as an application of our general abstract results, we present the solution to the mean-field LQ problem (which is a major
motivation of the current work). The open-loop optimal control for the MF-LQ control problem is characterized by the solution of a Fredholm type integral equation of the second kind. The theoretical
results are also applied to a concrete mean-variance portfolio selection problem arising from mathematical finance.

\ms

The rest of this paper is organized as follows. Some motivations of Problem (OLQ) are carefully presented in Section 2. Some general results for FSDEs and BSDEs with operator coefficients will be established as well. Section 3 is concerned with Problem (OLQ). Open-loop optimal controls are characterized, and the solvability of the relevant coupled FBSDEs with operator coefficients is established by the method of continuation. In Section 4, an MF-LQ optimal control problem and a mean-variance portfolio selection problem are worked out. Finally, we wind up this paper in Section 5.

\ms

\section{Preliminaries}

\subsection{Motivations}

In this subsection, we look at some motivations of our Problem (OLQ). First of all, for the state equation \rf{state}, let us look at some special cases.

\ms

$\bullet$ The classical linear SDE:
$$\left\{\2n\ba{ll}
\ns\ds\cA(s)\xi=A(s)\xi,\q\cC(s)\xi=C(s)\xi,\qq\forall\xi\in L^2_{\cF_s}(\O;\dbR^n),\\
\ns\ds\cB(s)\eta=B(s)\eta,\q\cD(s)\eta=D(s)\eta,\qq\forall\eta\in
L^2_{\cF_s}(\O;\dbR^m),\ea\right.$$
with $A(\cd),C(\cd),B(\cd),D(\cd)$ being some matrix-valued processes.

\ms

$\bullet$ The case of simple mean-field SDE (MF-SDE):
$$\left\{\2n\ba{ll}
\ns\ds\cA(s)\xi=A(s)\xi+\bar A(s)\dbE[\wt A(s)\xi],\q\cC(s)\xi=C(s)\xi+\bar C(s)\dbE[\wt C(s)\xi],\q\forall\xi\in L^2_{\cF_s}(\O;\dbR^n),\\
\ns\ds\cB(s)\eta=B(s)\eta+\bar B(s)\dbE[\wt B(s)\eta],\q\cD(s)\eta=D(s)\eta+\bar D(s)\dbE[\wt D(s)\eta],\q\forall\eta\in L^2_{\cF_s}(\O;\dbR^m),\ea\right.$$
for some matrix-valued processes $A(\cd),\bar A(\cd),\wt A(\cd)$, etc. In the case that all the coefficients are deterministic (as in \cite{Yong 2013,Huang-Li-Yong 2015,Li-Sun-Yong 2016}), the above will be reduced to the following simpler form:
$$\left\{\2n\ba{ll}
\ns\ds\cA(s)\xi=A(s)\xi+\bar A(s)\dbE[\xi],\q\cC(s)\xi=C(s)\xi+\bar C(s)\dbE[\xi],\q\forall\xi\in L^2_{\cF_s}(\O;\dbR^n),\\
\ns\ds\cB(s)\eta=B(s)\eta+\bar B(s)\dbE[\eta],\q\cD(s)\eta=D(s)\eta+\bar D(s)\dbE[\eta],\q\forall\eta\in L^2_{\cF_s}(\O;\dbR^m),\ea\right.$$
for some matrix-valued deterministic functions $A(\cd),\bar A(\cd)$, etc.

\ms

$\bullet$ The extended MF-SDE:
\bel{intA}\left\{\2n\2n\ba{ll}
\ns\ds\ba{ll}
\ns\ds\cA(s)\xi=A(s)\xi+\int_\dbR\bar A_\l(s)\dbE\big[\wt A_\l(s)\xi\big]\m(d\l),\\
\ns\ds\cC(s)\xi=C(s)\xi+\int_\dbR\bar C_\l(s)\dbE\big[\wt C_\l(s)\xi\big]\m(d\l),\ea\qq\forall\xi\in L^2_{\cF_s}(\O;\dbR^n),\\
\ns\ds\ba{ll}
\ns\ds\cB(s)\eta=B(s)\eta+\int_\dbR\bar B_\l(s)\dbE\big[\wt B_\l(s)\eta\big]\m(d\l),\\
\ns\ds\cD(s)\eta=D(s)\eta+\int_\dbR\bar D_\l(s)\dbE\big[\wt D_\l(s)\eta\big]\m(d\l),\ea\qq\forall\eta\in L^2_{\cF_s}(\O;\dbR^m),\ea\right.\ee
with $A(\cd),C(\cd),B(\cd),D(\cd)$ being matrix-valued processes and $\bar A_\l(\cd),\wt A_\l(\cd)$, etc. being families of matrix-valued processes parameterized by $\l\in\dbR$, and $\m(\cd)$ being a Borel measure on $\dbR$. Some conditions are needed in order the above make sense. A special case of the above is the following (with $\m(\cd)$ supported at $\{1,2,3,\cds\}$):
\bel{sumA}\left\{\2n\2n\ba{ll}
\ns\ds\ba{ll}
\ns\ds\cA(s)\xi=A(s)\xi+\sum_{k\ges1}\bar A_k(s)\dbE\big[\wt A_k(s)\xi\big]\equiv A(s)\xi+\bar\BA(s)^\top\dbE\big[\wt\BA(s)\xi\big],\\
\ns\ds\cC(s)\xi=C(s)\xi+\sum_{k\ges1}\bar C_k(s)\dbE\big[\wt C_k(s)\xi\big]\equiv C(s)\xi+\bar\BC(s)^\top\dbE\big[\wt\BC(s)\xi\big],\ea\qq\forall\xi\in L^2_{\cF_s}(\O;\dbR^n),\\
\ns\ds\ba{ll}
\ns\ds\cB(s)\eta=B(s)\eta+\sum_{k\ges1}\bar B_k(s)\dbE\big[\wt B_k(s)\eta\big]\equiv B(s)\eta+\bar\BB(s)^\top\dbE\big[\wt\BB(s)\eta\big],\\
\ns\ds\cD(s)\eta=D(s)\eta+\sum_{k\ges1}\bar D_k(s)\dbE\big[\wt D_k(s)\eta\big]\equiv D(s)\eta+\bar\BD(s)^\top\dbE\big[\wt\BD(s)\eta\big],\ea\qq\forall\eta\in L^2_{\cF_s}(\O;\dbR^m),\ea\right.\ee
with
\bel{BA}\left\{\2n\ba{ll}
\ns\ds\bar\BA(s)^\top=(\bar A_1(s),\bar A_2(s),\cds),\q\wt\BA(s)^\top=(\wt A_1(s)^\top,\wt A_2(s)^\top,\cds),\\
\ns\ds\bar\BC(s)^\top=(\bar C_1(s),\bar C_2(s),\cds),\q\wt\BC(s)^\top=(\wt C_1(s)^\top,\wt C_2(s)^\top,\cds),\\
\ns\ds\bar\BB(s)^\top=(\bar B_1(s),\bar B_2(s),\cds),\q\wt\BB(s)^\top=(\wt B_1(s)^\top,\wt B_2(s)^\top,\cds),\\
\ns\ds\bar\BD(s)^\top=(\bar D_1(s),\bar D_2(s),\cds),\q\wt\BD(s)^\top=(\wt D_1(s)^\top,\wt D_2(s)^\top,\cds),\ea\right.\ee
for some matrix-valued processes $\bar A_k(\cd),\wt A_k(\cd),\bar C_k(\cd),\wt C_k(\cd),\bar B_k(\cd),\wt B_k(\cd),\bar D_k(\cd),\wt D_k(\cd)$, $k=1,2,\cds$. We will look at the above case in details in Section 4.

\ms

From the above, we see that by allowing
$\cA(\cd),\cC(\cd),\cB(\cd),\cD(\cd)$ to be operator-valued
processes (not just matrix-valued processes), our state equation can
cover a big class of stochastic linear systems.

\ms

Next, we look at the cost functional. To get some feeling about the operators in the cost functional, let us look at the case compatible with \rf{sumA}--\rf{BA}. Let
$$\wt\BX(T)=\begin{pmatrix}X(T)\\ \dbE[\wt\BG X(T)]\end{pmatrix},\qq
\BX(s)=\begin{pmatrix}X(s)\\ \dbE[\wt\BQ(s)X(s)]\end{pmatrix},\qq
\Bu(s)=\begin{pmatrix}u(s)\\ \dbE[\wt\BR(s)u(s)]\end{pmatrix},$$
for some
\bel{BG}\wt\BG=\begin{pmatrix}\wt G_1\\ \wt G_2\\ \vdots\end{pmatrix},\qq
\wt\BQ(s)=\begin{pmatrix}\wt Q_1(s)\\ \wt Q_2(s)\\ \vdots\end{pmatrix},\qq
\wt\BR(s)=\begin{pmatrix}\wt R_1(s)\\ \wt R_2(s)\\ \vdots\end{pmatrix}.\ee
For the terminal cost, we propose the quadratic term as
$$\ba{ll}
\ns\ds\dbE\lan\BG\wt\BX(T),\wt\BX(T)\ran=\dbE\bigg\lan\begin{pmatrix}G&\bar\BG^\top\\ \bar\BG&\h\BG\end{pmatrix}\begin{pmatrix}X(T)\\ \dbE[\wt\BG X(T)]\end{pmatrix},\begin{pmatrix}X(T)\\ \dbE[\wt\BG X(T)]\end{pmatrix}\bigg\ran\\
\ns\ds=\dbE\[\lan GX(T),X(T)\ran+\lan\bar\BG^\top\dbE[\wt\BG X(T)],X(T)\ran+\lan\bar\BG X(T),\dbE[\wt\BG X(T)]\ran+\lan\h\BG\dbE[\wt\BG X(T)],
\dbE[\wt\BG X(T)]\ran\]\\
%
%
%
\ns\ds=\dbE\blan GX(T)+\sum_{k\ges1}\(\bar G_k^\top\dbE[\wt G_kX(T)]+\wt G_k^\top\dbE[\bar G_kX(T)]\)+\sum_{i,j\ges1}\wt G_j^\top\dbE[\h G_{ij}]\dbE[\wt G_iX(T)],X(T)\bran\\
\ns\ds\equiv\dbE\lan GX(T)+\bar\BG^\top\dbE[\wt\BG X(T)]+\wt\BG^\top\dbE[\bar\BG X(T)]
+\wt\BG^\top\dbE[\h\BG]\dbE[\wt\BG X(T)],X(T)\ran\equiv\dbE\lan\cG X(T),X(T)\ran,\ea$$
and the linear term as
$$\ba{ll}
\ns\ds\dbE\lan\Bg,\wt\BX(T)\ran=\dbE\bigg\lan\begin{pmatrix}g_0\\ \bar\Bg\end{pmatrix},\begin{pmatrix}X(T)\\ \dbE[\wt\BG X(T)]\end{pmatrix} \bigg\ran
=\dbE\[\lan g_0,X(T)\ran+\lan\bar\Bg,\dbE[\wt\BG X(T)]\ran\]\\
\ns\ds=\dbE\[\lan g_0,X(T)\ran+\sum_{k\ges1}\lan\bar g_k,\dbE[\wt G_kX(T)]\ran\]=\dbE\lan g_0+\sum_{k\ges1}\wt G_k^\top\dbE[\bar g_k],X(T)\ran\\
\ns\ds=\dbE\lan g_0+\wt\BG^\top\dbE[\Bg],X(T)\ran\equiv\dbE\lan g,X(T)\ran,\ea$$
for some
\bel{G}G^\top=G,\qq\bar\BG=\begin{pmatrix}\bar G_1\\ \bar G_2\\ \vdots\end{pmatrix},\qq
\h\BG=\begin{pmatrix}\h G_{11}&\h G_{12}&\cds\\
                        \h G_{21}&\h G_{22}&\cds\\
                        \vdots&\vdots&\ddots\end{pmatrix},\qq\h\BG^\top=\h\BG,\qq\bar\Bg=\begin{pmatrix}\bar g_1\\ \bar g_2\\ \vdots\end{pmatrix}.\ee
For the running cost, we propose the quadratic terms as
$$\ba{ll}
\ns\ds\dbE\lan\BQ(s)\BX(s),\BX(s)\ran=\dbE\bigg\lan\begin{pmatrix}Q(s)&\bar\BQ(s)^\top\\ \bar\BQ(s)&\h\BQ(s)\end{pmatrix}\begin{pmatrix}X(s)\\ \dbE[\wt{\BQ}(s)X(s)]\end{pmatrix},\begin{pmatrix}X(s)\\ \dbE[\wt{\BQ}(s)X(s)]\end{pmatrix}\bigg\ran\\
\ns\ds=\dbE\[\lan Q(s)X(s),X(s)\ran+\lan\bar\BQ(s)^\top \dbE[\wt\BQ(s)X(s)],X(s)\ran+\lan\bar\BQ(s)X(s),\dbE[\wt\BQ(s)X(s)]\ran\\
\ns\ds\qq\qq\qq\qq\qq\qq\qq+\lan\h\BQ(s)\dbE[\wt\BQ(s)X(s)],
\dbE[\wt\BQ(s)X(s)]\ran\]\\
%
%
%
\ns\ds=\dbE\[\lan Q(s)X(s)+\sum_{k\ges1}\(\bar Q_k(s)^\top\dbE[\wt Q_k(s)X(s)]+
\wt Q_k(s)^\top\dbE[\bar Q_k(s)X(s)]\)\\
\ns\ds\qq\qq\qq\qq\qq\qq\qq+\sum_{i,j\ges1}\wt Q_j(s)^\top\dbE[\h Q_{ij}(s)]\dbE[\wt Q_i(s)X(s)],X(s)\ran\]\\
\ns\ds=\dbE\lan Q(s)X(s)+\bar\BQ(s)^\top\dbE[\wt\BQ(s)X(s)]+\wt\BQ(s)^\top\dbE[\bar\BQ(s)X(s)]
+\wt\BQ(s)^\top\dbE[\h\BQ(s)]\dbE[\wt\BQ(s)X(s)],X(s)\ran\\
\ns\ds\equiv\dbE\lan\cQ(s)X(s),X(s)\ran,\ea$$
and similarly,
$$\ba{ll}
\ns\ds\dbE\lan\BR(s)\Bu(s),\Bu(s)\ran=\dbE\bigg\lan\begin{pmatrix}R(s)&\bar\BR(s)^\top\\ \bar\BR(s)&\h\BR(s)\end{pmatrix}\begin{pmatrix}u(s)\\ \dbE[\wt\BR(s)u(s)]\end{pmatrix},\begin{pmatrix}u(s)\\ \dbE[\wt\BR(s)u(s)]\end{pmatrix}\bigg\ran\\
\ns\ds=\dbE\[\lan R(s)u(s),u(s)\ran+\lan\bar\BR(s)^\top\dbE[\wt\BR(s)u(s)],u(s)\ran+\lan\bar\BR(s)u(s),\dbE[\wt\BR(s)u(s)]\ran\\
\ns\ds\qq\qq\qq\qq\qq\qq+\lan\h\BR(s)\dbE[\wt\BR(s)u(s)],\dbE[\wt\BR(s)u(s)]\ran\]\\
\ns\ds=\dbE\[\lan R(s)u(s)+\sum_{k\ges1}\(\lan\bar R_k(s)^\top\dbE[\wt R_k(s)u(s)]+\wt R_k(s)\dbE[\bar R_k(s)u(s)]\)\\
\ns\ds\qq\qq\qq\qq\qq\qq  +\sum_{i,j\ges1}\wt R_j(s)^\top\dbE[\h R_{ij}(s)]\dbE[\wt R_i(s)u(s)],u(s)\ran\]\\
\ns\ds=\dbE\lan R(s)u(s)+\bar\BR(s)^\top\dbE[\wt\BR(s)u(s)]+\wt\BR(s)^\top\dbE[\bar\BR(s)u(s)]
+\wt\BR(s)^\top\dbE[\h\BR(s)]\dbE[\wt\BR(s)u(s)],u(s)\ran\\
\ns\ds\equiv\dbE\lan\cR(s)u(s),u(s)\ran,\ea$$
$$\ba{ll}
\ns\ds\dbE\lan\BS(s)\BX(s),\Bu(s)\ran=\dbE\bigg\lan\begin{pmatrix}S(s)&\wt\BS(s)^\top\\ \bar\BS(s)&\h\BS(s)\end{pmatrix}\begin{pmatrix}X(s)\\ \dbE[\wt\BQ(s)X(s)]\end{pmatrix},\begin{pmatrix}u(s)\\ \dbE[\wt\BR(s)u(s)]\end{pmatrix}\bigg\ran\\
\ns\ds=\dbE\[\lan S(s)X(s),u(s)\ran+\lan\wt\BS(s)^\top\dbE[\wt\BQ(s)X(s)],u(s)\ran+\lan\bar\BS(s)X(s),\dbE[\wt\BR(s)u(s)]\ran\\
\ns\ds\qq\qq\qq\qq\qq\qq+\lan\h\BS(s)\dbE[\wt\BQ(s)X(s)],\dbE[\wt\BR(s)u(s)]\ran\]\\
\ns\ds=\dbE\[\lan S(s)X(s)+\sum_{k\ges1}\(\wt \BS_k(s)^\top\dbE[\wt Q_k(s)X(s)]+\wt R_k(s)\dbE[\bar S_k(s)X(s)]\)\\
\ns\ds\qq\qq\qq\qq\qq\qq +\sum_{i,j\ges1}\wt R_i(s)^\top\dbE[\h S_{ij}(s)]\dbE[\wt Q_(s)X(s)],u(s)\ran\]       \\
\ns\ds=\dbE\lan S(s)X(s)+\wt\BS(s)^\top\dbE[\wt\BQ(s)X(s)]+\wt\BR(s)^\top\dbE[\bar\BS(s)X(s)]
+\wt\BR(s)^\top\dbE[\h\BS(s)]\dbE[\wt\BQ(s)X(s)],u(s)\ran\\
\ns\ds\equiv\dbE\lan\cS(s)X(s),u(s)\ran,\ea$$
and propose the linear terms as
$$\ba{ll}
\ns\ds\dbE\lan\Bq(s),\BX(s)\ran=\dbE\bigg\lan\begin{pmatrix}q_0(s)\\ \bar\Bq(s)\end{pmatrix},
\begin{pmatrix}X(s)\\ \dbE[\wt\BQ(s)X(s)]\end{pmatrix}\bigg\ran=\dbE\[\lan q_0(s),X(s)\ran+\lan\bar\Bq(s),\dbE[\wt\BQ(s)X(s)]\ran\]\\
\ns\ds=\dbE\[\lan q_0(s),X(s)\ran+\sum_{k\ges1}\lan\bar q_k(s),\dbE\wt Q_k(s)X(s)]\ran\]=\dbE\lan q_0(s)+\sum_{k\ges1}\wt Q_k(s)^\top\dbE[\bar q_k(s)],X(s)\ran\\
\ns\ds=\dbE\lan q_0(s)+\wt\BQ(s)^\top\dbE[\bar\Bq(s)],X(s)\ran\equiv\dbE\lan q(s),X(s)\ran,\\
\ns\ds\dbE\lan\Brho(s),\Bu(s)\ran=\dbE\bigg\lan\begin{pmatrix}\rho_0(s)\\ \bar\Brho(s)\end{pmatrix},
\begin{pmatrix}u(s)\\ \dbE[\wt\BR(s)u(s)]\end{pmatrix}\bigg\ran=\dbE\[\lan\rho_0(s),u(s)\ran+\lan\bar\Brho(s),
\dbE[\wt\BR(s)u(s)]\ran\]\\
\ns\ds=\dbE\[\lan\rho_0(s),X(s)\ran+\sum_{k\ges1}\lan\bar\rho_k(s),\dbE[\wt R_k(s)u(s)]\]=\dbE\lan\rho_0(s)+\sum_{k\ges1}\wt R_k(s)^\top\dbE[\bar\rho_k(s)],u(s)\ran\\
\ns\ds=\dbE\lan \rho_0(s)+\wt\BR(s)^\top\dbE[\bar\Brho(s)],u(s)\ran\equiv\dbE\lan\rho(s),u(s)\ran,\ea$$
for some
\bel{Q}Q(s)^\top=Q(s),\q\bar\BQ(s)=\begin{pmatrix}\bar Q_1(s)\\ \bar Q_2(s)\\ \vdots\end{pmatrix},\q
\h\BQ(s)=\begin{pmatrix}\h Q_{11}(s)&\h Q_{12}(s)&\cds\\
                        \h Q_{21}(s)&\h Q_{22}(s)&\cds\\
                        \vdots&\vdots&\ddots\end{pmatrix},\q\h\BQ(s)^\top=\h\BQ(s),\ee
\bel{R}R(s)^\top=R(s),\q\bar\BR(s)=\begin{pmatrix}\bar R_1(s)\\ \bar R_2(s)\\ \vdots\end{pmatrix},\q
\h\BR(s)=\begin{pmatrix}\h R_{11}(s)&\h R_{12}(s)&\cds\\
                        \h R_{21}(s)&\h R_{22}(s)&\cds\\
                        \vdots&\vdots&\ddots\end{pmatrix},\q\h\BR(s)^\top=\h\BR(s).\ee
\bel{S}\bar\BS(s)=\begin{pmatrix}\bar S_1(s)\\ \bar S_2(s)\\ \vdots\end{pmatrix},\qq
\wt\BS(s)=\begin{pmatrix}\wt S_1(s)\\ \wt S_2(s)\\ \vdots\end{pmatrix},\qq
\h\BS(s)=\begin{pmatrix}\h S_{11}(s)&\h S_{12}(s)&\cds\\
                        \h S_{21}(s)&\h S_{22}(s)&\cds\\
                        \vdots&\vdots&\ddots\end{pmatrix}.\ee
\bel{g}\qq
\bar\Bq(s)=\begin{pmatrix}\bar q_1(s)\\ \bar q_2(s)\\ \vdots\end{pmatrix},\qq
\bar\Brho(s)=\begin{pmatrix}\bar\rho_1(s)\\ \bar\rho_2(s)\\ \vdots\end{pmatrix}.\ee
We see that, in the above case,
\bel{Sec2.1_Cost_Operators}
\left\{\2n\ba{ll}
\ns\ds
\cG\xi=G\xi+\wt\BG^\top\dbE[\bar\BG\xi]+\bar\BG^\top\dbE[\wt\BG\xi]
+\wt\BG^\top\dbE[\h\BG]\dbE[\wt\BG\xi],\q\xi\in L^2(\O;\dbR^n),\\
\ns\ds
\cQ(s)\xi\1n=\1n Q(s)\xi\1n+\1n\wt\BQ(s)^\top\dbE[\bar\BQ(s)\xi]\1n+\1n\bar\BQ(s)^\top\dbE[\wt\BQ(s)\xi]
\1n+\1n\wt\BQ(s)^\top\dbE[\h\BQ(s)]\dbE[\wt\BQ(s)\xi],\q\xi\in
L^2_{\cF_s}(\O;\dbR^n),\\
\ns\ds
\cS(s)\xi=S(s)\xi+\wt\BR(s)^\top\dbE[\bar\BS(s)\xi]+\wt\BS(s)^\top\dbE[\wt\BQ(s)\xi]
+\wt\BR(s)^\top\dbE[\h\BS(s)]\dbE[\wt\BQ(s)\xi],\q\xi\in
L^2_{\cF_s}(\O;\dbR^n),\\
\ns\ds\cR(s)\eta\1n=\1n R(s)\eta\1n+\1n\wt\BR(s)^\top\dbE[\bar\BR(s)\eta]\1n+\1n\bar\BR(s)^\top\dbE[\wt\BR(s)\eta]
\1n+\1n\wt\BR(s)^\top\dbE[\h\BR(s)]\dbE[\wt\BR(s)\eta],\q\eta\in
L^2_{\cF_s}(\O;\dbR^m),\\
\ns\ds  g=g_0+\wt\BG^\top\dbE[\bar\Bg],\q
q(s)=q_0(s)+\wt\BQ(s)^\top\dbE[\bar\Bq(s)],\q
\rho(s)=\rho_0(s)+\wt\BR(s)^\top\dbE[\bar\Brho(s)].\ea\right.\ee
When all the weighting functions in the cost functional are deterministic, the above will be reduced to the following: (see \cite{Yong 2013})
$$\left\{\2n\ba{ll}
\ns\ds\cG\xi=G\xi+\bar G\dbE[\xi],\q\xi\in L^2(\O;\dbR^n),\\
\ns\ds\cQ(s)\xi=Q(s)\xi+\bar Q(s)\dbE[\xi],\q\cS(s)\xi=S(s)\xi+\bar S(s)\dbE[\xi],\qq\xi\in L^2_{\cF_s}(\O;\dbR^n),\\
\ns\ds\cR(s)\eta=R(s)\eta+\bar R(s)\dbE[\eta],\q\eta\in L^2_{\cF_s}(\O;\dbR^m).\ea\right.$$
The above suggests that if the coefficients of the state equation are given by \rf{intA}, the corresponding operators in the cost functional could look like the following:
\bel{}\left\{\2n\ba{ll}
\ns\ds\cG\xi=G\xi+\2n\int\(\wt G_\l^\top\dbE[\bar G_\l\xi]+\bar G^\top_\l\dbE[\wt G_\l\xi]\)\m(d\l)
+\2n\int\3n\int\wt G_\l^\top\h G_{\l\n}\dbE[\wt G_\n\xi]\m(d\l)\m(d\n),\q\xi\in L^2(\O;\dbR^n),\\
\ns\ds\cQ(s)\xi=Q(s)\xi+\int\(\wt Q_\l(s)^\top\dbE[\bar Q_\l(s)\xi]+\bar Q_\l(s)^\top
\dbE[\wt Q_\l(s)\xi]\)\m(d\l)\\
\ns\ds\qq\qq\qq\qq+\int\3n\int\wt Q_\l(s)^\top\h Q_{\l\n}(s)\dbE[\wt Q_\n(s)\xi]\m(d\l)\m(d\n),\qq
\xi\in L^2_{\cF_s}(\O;\dbR^n),\\
\ns\ds\cS(s)\xi=S(s)\xi+\int\(\wt R_\l(s)^\top\dbE[\bar S_\l(s)\xi]+\wt S_\l(s)^\top\dbE[\wt Q(s)\xi]\)\m(d\l)\\
\ns\ds\qq\qq\qq\qq+\int\3n\int\wt R_\l(s)^\top\h S_{\l\n}(s)\dbE[\wt Q_\n(s)\xi]\m(d\l)\m(d\n),\qq\xi\in L^2_{\cF_s}(\O;\dbR^n),\\
\ns\ds\cR(s)\eta=R(s)\eta+\int\(\wt R_\l(s)^\top\dbE[\bar R_\l(s)\eta]+\bar R_\l(s)^\top
\dbE[\wt R_\l(s)\eta]\)\m(d\l)\\
\ns\ds\qq\qq\qq\qq+\int\3n\int\wt R_\l(s)^\top\h R_{\l\n}(s)\dbE[\wt R_\n(s)\eta]\m(d\l)\m(d\n),\qq\eta\in L^2_{\cF_s}(\O;\dbR^m).\ea\right.\ee
In the above, $\h G_{\l\n},\h Q_{\l\n}(\cd),\h S_{\l\n}(\cd),\h R_{\l\n}(\cd)$ are deterministic, and
\bel{}\h G_{\l\n}^\top=\h G_{\n\l},\qq\h Q_{\l\n}(s)^\top=\h Q_{\n\l}(s),\qq\h R_{\l\n}(s)^\top=\h R_{\n\l}(s),\qq\forall\l,\n\in\dbR,~s\in[0,T].\ee

The above shows that our framework can cover many problems involving
mean fields.

\subsection{The state equation and the cost functional}

We now return to our state equation \rf{state} and cost functional \rf{cost}. Recall the Hilbert space
$$L^2_{\cF_s}(\O;\dbR^n)\subseteq L^2(\O;\dbR^n)\equiv L^2_{\cF_T}(\O;\dbR^n),\qq s\in[0,T],$$
with the norm:
$$\|\xi\|_2=\(\dbE|\xi|^2\)^{1\over2},\qq\forall\xi\in L^2_{\cF_s}(\O;\dbR^n).$$
Also, we recall the spaces $L^2_\dbF(\O;L^1(t,T;\dbR^n))$ and
$L^2_\dbF(\O;C([t,T];\dbR^n))$ introduced in the first section. We
now introduce the following definition concerning the
operator-valued processes.

\bde{} \rm An operator-valued process  $\cB:[0,T]\to\sL\(L^2(\O;\dbR^m);L^2(\O;\dbR^n)\)$ is said to be {\it strongly $\dbF$-progressively measurable} if
$$\cB(s)\in\sL\(L^2_{\cF_s}(\O;\dbR^m);L^2_{\cF_s}(\O;\dbR^n)\),\qq\forall s\in[0,T],$$
and for any $\eta(\cd)\in L^2_\dbF(0,T;\dbR^m)$, $\cB(\cd)\eta(\cd)$ is $\dbF$-progressively measurable. The set of all strongly $\dbF$-progressively measurable operator-valued processes valued in
$\sL\(L^2(\O;\dbR^m);L^2(\O;\dbR^n)\)$ is denoted by $\sL_\dbF\(L^2(\O;\dbR^m)$; $L^2(\O;\dbR^n)\)$, and
denote
$$\sL_\dbF\(L^2(\O;\dbR^n);L^2(\O;\dbR^n)\)=\sL_\dbF\(L^2(\O;\dbR^n)\).$$
Further, for any $p\in[1,\i]$, we denote
$$\sL^p_\dbF\(L^2(\O;\dbR^m);L^2(\O;\dbR^n)\)=\Big\{\cB(\cd)\in\sL_\dbF\(L^2(\O;\dbR^m);L^2(\O;\dbR^n)\)
\bigm|\|\cB(\cd)\|_p<\i\Big\},$$
where
$$\|\cB(\cd)\|_p=\left\{\2n\ba{ll}
\ns\ds\(\int_0^T\|\cB(t)\|^pdt\)^{1\over p},\qq\qq p\in[1,\i),\\
\ns\ds\esssup_{t\in[0,T]}\|\cB(t)\|, \qq\qq\qq p=\i,\ea\right.$$
with $\|\cB(t)\|$ being the operator norm of $\cB(t)$ for given $t\in[0,T]$, defined by
$$\|\cB(t)\|=\sup\bigg\{\(\dbE\big[|\cB(t)\eta|^2\big]\)^{1\over2}\Bigm|\eta\in L^2_{\cF_t}(\O;\dbR^m),\ \(\dbE\big[|\eta|^2\big]\)^{1\over2}=1\bigg\}.$$
Also, we denote
$$\sL^p_\dbF\(L^2(\O;\dbR^n);L^2(\O;\dbR^n)\)=\sL_\dbF^p\(L^2(\O;\dbR^n)\),$$
and
$$\sS^p_\dbF\Big(L^2(\O;\dbR^n)\Big)=\Big\{\cB(\cd)\in\sL^p_\dbF\Big(L^2(\O;\dbR^n)\Big)\bigm|\cB(s)\in \sS\Big(L^2_{\cF_s}(\O;\dbR^n) \Big),\ \forall s\in [0,T]\Big\}.$$

\ede

Strong measurability for operator-valued functions can be found in \cite{Yosida 1980}.
Our operator-valued processes have an additional feature of $\dbF$-progressive measurability.
Therefore, the above definition is necessary.

\ms

Now, for state equation \rf{state}, we introduce the following hypothesis.

\ms

{\bf(H1)} Let
$$\ba{ll}
\ns\ds\cA(\cd)\in\sL^1_\dbF\(L^2(\O;\dbR^n)\),\q\cB(\cd)\in\sL^2_\dbF\(L^2(\O;\dbR^m);L^2(\O;\dbR^n)\),\\
\ns\ds\cC(\cd)\in\sL^2_\dbF\(L^2(\O;\dbR^n)\),\q\cD(\cd)\in\sL^\i_\dbF\(L^2(\O;\dbR^m);L^2(\O;\dbR^n)\),\\
\ns\ds b(\cd)\in L^2_\dbF(\O;L^1(0,T;\dbR^n)),\q\si(\cd)\in
L^2_\dbF(0,T;\dbR^n).\ea$$

\ms

The following result gives the well-posedness of \rf{state}.

\bp{Well-posedness-state} \sl Let {\rm(H1)} hold. Then for any $(t,x)\in\sD$, and $u(\cd)\in\cU[t,T]$, there exists a unique solution $X(\cd)\equiv X(\cd\,;t,x)\in L^2_\dbF(\O;C([t,T];\dbR^n))$ to \rf{state} and the following estimate holds:
\bel{Yu_Sec2.2_FSDE_Estimate}\ba{ll}
\ns\ds\dbE\[\sup_{s\in [t,T]}|X(s)|^2\]\les K\Big[\|x\|_2^2+\(\|\cB(\cd)\|_2^2+\|\cD(\cd)\|_\i^2\)\int_t^T\|u(r)\|_2^2dr\\
\ns\ds\qq\qq\qq\qq\qq\qq\qq+\Big\| \int_t^T |b(r)|dr
\Big\|_2^2+\int_t^T\|\si(r)\|_2^2dr\],\ea\ee
for some constant $K>0$ depending on $\|\cA(\cd)\|_1$ and $\|\cC(\cd)\|_2$.

\ep

\it Proof. \rm Let $(t,x)\in\sD$ be fixed and $u(\cd)\in\cU[t,T]$ be given. For any $\bar X(\cd)\in L^2_\dbF(\O;C([t,T];\dbR^n))$, we define the process $X(\cd)$ by
$$X(s)=x+\int_t^s\Big(\cA(r)\bar X(r)+\cB(r)u(r)+b(r)\Big)dr+\int_t^s\Big(\cC(r)\bar X(r)+\cD(r)u(r)+\si(r)\Big)dW(r),\q
s\in [t,T].$$
For any $t\les t_1<t_2\les T$ and any $s\in[t_1,t_2]$,
$$\ba{ll}
\ns\ds|X(s)|^2\les3\Big\{|X(t_1)|^2+\(\int_{t_1}^{s}\big|\cA(r)\bar X(r)+\cB(r)u(r)+b(r)\big|dr\)^2\\
\ns\ds\qq\qq\qq\qq+\Big|
\int_{t_1}^s\big(\cC(r)\bar X(r)+\cD(r)u(r)+\si(r)\big)dW(r)\Big|^2\Big\}.\ea$$
By Minkowski's integral inequality and the Burkholder--Davis--Gundy
inequality (see Karatzas-Shreve \cite{Karatzas-Shreve
1991}, Chapter 3, Theorem 3.28), we have
$$\ba{ll}
\ns\ds\dbE\[\sup_{\t\in [t_1,s]}|X(\t)|^2\]\les3\Big\{\|X(t_1)\|_2^2+\Big\|\int_{t_1}^s \big|\cA(r)\bar X(r)+\cB(r)u(r)+b(r)\big|dr\Big\|_2^2\\
\ns\ds\qq\qq\qq\qq\qq\qq\qq+\dbE\[\sup_{\t\in[t_1,s]}\Big|\int_{t_1}^\t\(\cC(r)\bar X(r)+\cD(r)u(r)+\si(r)\Big)dW(r)\Big|^2 \] \Big\}\\
\ns\ds \les3\Big\{\|X(t_1)\|_2^2+3\( \Big\|\int_{t_1}^s |\cA(r)\bar X(r)|dr\Big\|_2^2 +\Big\|\int_{t_1}^s |\cB(r)u(r)|dr\Big\|_2^2 +\Big\|\int_{t_1}^s |b(r)| dr\Big\|_2^2 \)\\
\ns\ds\qq\qq\qq\qq\qq\qq\qq
+\dbE\[\sup_{\t\in[t_1,s]}\Big|\int_{t_1}^\t\(\cC(r)\bar
X(r)+\cD(r)u(r)+\si(r)\Big)dW(r)\Big|^2 \] \Big\}\\
\ns\ds\les3\Big\{\|X(t_1)\|_2^2+3\[ \(\int_{t_1}^s \|\cA(r)\bar X(r)\|_2 dr\)^2+\(\int_{t_1}^s\|\cB(r)u(r)\|_2dr\)^2 +\Big\|\int_{t_1}^s |b(r)|dr\Big\|_2^2 \]\\
\ns\ds\qq\qq\qq\qq\qq\qq\qq+3c_2\int_{t_1}^s\Big(\|\cC(r)\bar X(r)\|_2^2+\|\cD(r)u(r)\|_2^2+\|\si(r)\|_2^2\Big)dr\Big\}\\
\ns\ds\les3\Big\{\|X(t_1)\|_2^2+3\[\(\int_{t_1}^s\|\cA(r)\|dr\)^2+c_2\int_{t_1}^s\|\cC(r)\|^2ds\]
\sup_{r\in[t_1,s]}\|\bar X(r)\|_2^2\\
\ns\ds\q+3\(\int_{t_1}^s\|\cB(r)\|^2dr+c_2\sup_{r\in[t_1,s]}\|\cD(r)\|^2\)\int_{t_1}^s\|u(r)\|_2^2dr
+3\Big\| \int_{t_1}^s |b(r)|dr
\Big\|_2^2+3c_2\int_{t_1}^s\|\si(r)\|_2^2dr\Big\},\ea$$
%
where $c_2>0$ is the constant in the Burkholder--Davis--Gundy inequality. Thus, it follows that for some generic constant $K>0$, (hereafter, $K>0$ will stand for a generic constant which could be different from line to line)
$$\ba{ll}
\ns\ds\dbE\[\sup_{\t\in [t_1,s]}|X(\t)|^2\]\les K
\Big\{\|X(t_1)\|_2^2+\Big\|\int_{t_1}^s |b(r)| dr \Big\|_2^2+\int_{t_1}^s\|\si(r)\|_2^2dr\\
\ns\ds\qq\qq\qq\qq\qq+\(\int_{t_1}^s\|\cB(r)\|^2dr+\sup_{r\in[t_1,s]}\|\cD(r)\|^2\)\int_{t_1}^s
\|u(r)\|_2^2dr\Big\}\\
\ns\ds\qq\qq\qq\qq\qq+K\[\(\int_{t_1}^s\|\cA(r)\|dr\)^2+\int_{t_1}^s\|\cC(r)\|^2dr\] \dbE\[\sup_{\t\in[t_1,s]}|\bar X(\t)|^2\].\ea$$
For $\d>0$ small enough, we have
$$K\[\(\int_t^{t+\d}\|\cA(r)\|dr\)^2+\int_t^{t+\d}\|\cC(r)\|^2dr\]<1.$$
Thus, the map $\bar X(\cd)\mapsto X(\cd)$ is a contraction from $L^2_\dbF(\O;C([t,t+\d];\dbR^n))$ into itself. Therefore, it admits a unique fixed point which is a solution to the state equation on $[t,t+\d]$. Repeating the same argument, we can obtain the unique existence of the solution $X(\cd)\in
L^2_\dbF(\O;C([t,T];\dbR^n))$ to the state equation \rf{state}.

\ms

Moreover, for the solution $X(\cd)$, from the above, we have
$$\ba{ll}
\ns\ds\dbE\[\sup_{\t\in[t,s]}|X(\t)|^2\]\les K\Big\{|x|^2+\Big\|
\int_t^s |b(r)|dr \Big\|_2^2
+\int_t^s\|\si(r)\|_2^2dr+\(\|\cB(\cd)\|_2^2+\|\cD(\cd)\|_\i^2\)\int_{t_1}^s
\|u(r)\|_2^2dr\Big\}\\
\ns\ds\qq\qq\qq\qq\q+K\[\(\int_t^s\|\cA(r)\|dr\)^2+\int_t^s\|\cC(r)\|^2dr\]
\dbE\[\sup_{\t\in[t,r]}|X(\t)|^2\].\ea$$
Then by a simple iteration argument, we obtain \rf{Yu_Sec2.2_FSDE_Estimate}. \endpf

\ms

Now, for the cost functional, we introduce the following hypothesis.

\ms

{\bf(H2)} The operator $\cG\in\sS\(L^2(\O;\dbR^n)\)$ and the operator-valued processes
$$\cQ(\cd)\in\sS^1_\dbF\(L^2(\O;\dbR^n)\),\q
\cS(\cd)\in\sL^2_\dbF\(L^2(\O;\dbR^n);L^2(\O;\dbR^m)\),\q
\cR(\cd)\in\sS^\i_\dbF\(L^2(\O;\dbR^m)\).$$
Also,
$$g\in L^2(\O;\dbR^n),\q q(\cd)\in L^2_\dbF(\O;L^1(0,T;\dbR^n)),\q\rho(\cd)\in L^2_\dbF(0,T;\dbR^m).$$

We have the following result.

\bp{} \sl Let {\rm(H1)--(H2)} hold. Then for any $(t,x)\in\sD$ and any $u(\cd)\in\cU[t,T]$, the cost functional $J(t,x;u(\cd))$ is well-defined.

\ep

\it Proof. \rm First of all, Proposition \ref{Well-posedness-state} implies that the state
equation \rf{state} admits a unique state process $X(\cd)\equiv X(\cd;t,x,u(\cd))\in L^2_\dbF(\O;C([t,T];\dbR^n))$. Let us observe the following estimates:
$$\ba{ll}
\dbE|\lan\cG X(T),X(T)\ran|\les\|\cG\|\|X(T)\|_2^2,\q\dbE|\lan g,X(T)\ran|\les\|g\|_2\|X(T)\|_2;\\
\ns\ds\int_t^T\dbE|\lan\cR(s)u(s),u(s)\ran|ds\les\int_t^T\|\cR(s)\|\|u(s)\|_2^2ds\les\|\cR(\cd)\|_\i \int_t^T\|u(s)\|_2^2ds;\\
\ns\ds\int_t^T\dbE|\lan\cS(s)X(s),u(s)\ran|ds\les\(\int_t^T\|\cS(s)X(s)\|_2^2ds\)^{1\over2} \(\int_t^T\|u(s)\|_2^2ds\)^{1\over2};\\
\ns\ds\qq\qq\qq\qq\qq\q\les\|\cS(\cd)\|_2\(\dbE\[\sup_{s\in [t,T]}|X(s)|^2\]\)^{1\over2}\(\int_t^T\|u(s)\|_2^2ds\)^{1\over2};\\
\ns\ds\int_t^T\dbE|\lan\cQ(s)X(s),X(s)\ran|ds\les\int_t^T\|\cQ(s)\|\dbE|X(s)|^2ds\les\|\cQ(\cd)\|_1
\dbE\[\sup_{s\in[t,T]}|X(s)|^2\];\ea$$
$$\ba{ll}
\ns\ds\int_t^T\dbE|\lan q(s),X(s)\ran|ds\les \dbE \[ \( \int_t^T |q(s)| ds \) \sup_{s\in [t,T]} |X(s)| \] \les  \Big\| \int_t^T |q(s)|ds \Big\|_2 \( \dbE\[ \sup_{s\in [t,T]}|X(s)|^2 \] \)^{\frac 1 2};\\
\ns\ds\int_t^T\dbE|\lan\rho(s),u(s)\ran|ds\les\(\int_t^T\|\rho(s)\|_2^2ds\)^{1\over2}
\(\int_t^T\|u(s)\|_2^2ds\)^{1\over2}.\ea$$
This implies that the cost functional $J(t,x;u(\cd))$ is well-defined. \endpf

\ms

Next we look at the quadratic form of the cost functional \rf{cost}, from which we will get some abstract results for Problem (OLQ).

\ms

It is clear that for given $t\in[0,T)$, $(x,u(\cd))\mapsto X(\cd\,;t,x,u(\cd))$ is affine. Therefore, we may write
$$X(\cd;t,x,u(\cd))=[F_1(t)u(\cd)](\cd)+[F_0(t)x](\cd)+f_0(t,\cd),$$
where
$$F_1(t)\in\sL\(\cU[t,T];L^2_\dbF(t,T;\dbR^n)\),\q F_0(t)\in\sL\(L^2_{\cF_t}(\O;\dbR^n);L^2_\dbF(t,T;\dbR^n)\),\q
f_0(t,\cd)\in L^2_\dbF(t,T;\dbR^n).$$
Let
$$\h F_1(t)u(\cd)=[F_1(t)u(\cd)](T),\q\h F_0(t)x=[F_0(t)x](T),\q\h f_0(t)=f_0(t,T).$$
Consequently,
\bel{J}\ba{ll}
\ns\ds J(t,x;u(\cd))\\
\ns\ds \1n=\1n\dbE\[\big\lan\cG\big\{\h F_1(t)u(\cd)\1n+\1n\h F_0(t)x\1n+\1n\h f_0(t)\big\},\h F_1(t)u(\cd)\1n+\1n\h F_0(t)x\1n+\1n\h f_0(t)\big\ran +2\big\lan g,\h F_1(t)u(\cd)+\h F_0(t)x+\h f_0(t)\big\ran\\
\ns\ds\q+\2n\int_t^T\3n\(\big\lan\cQ(s)\big\{[F_1(t)u(\cd)](s)+[F_0(t)x](s)+f_0(t,s)\big\},
[F_1(t)u(\cd)](s)+[F_0(t)x](s)+f_0(t,s)\big\ran\\
\ns\ds\qq\qq\qq\q+2\lan\cS(s)\big\{[F_1(t)u(\cd)](s)+[F_0(t)x](s)+f_0(t,s)\big\},u(s)\ran+\lan\cR(s)u(s),
u(s)\ran\\
\ns\ds\qq\qq\qq\q+2\big\lan q(s),[F_1(t)u(\cd)](s)+[F_0(t)x](s)+f_0(t,s)\ran+2\lan\rho(s),u(s)\big\ran\)ds\]\\
\ns\ds\equiv\big\lan\cG[\h F_1u+\h F_0x+\h f_0],\h F_1u+\h F_0x+\h f_0\big\ran+2\big\lan g,\h F_1u+\h F_0x+\h f_0\big\ran\\
\ns\ds\q+\big\lan\cQ(F_1u+F_0x+f_0),F_1u+F_0x+f_0\big\ran+2\big\lan\cS(F_1u+F_0x+f_0),u\big\ran+\big\lan\cR u,u\big\ran\\
\ns\ds\q+2\big\lan q,F_1u+F_0x+f_0\big\ran+2\big\lan\rho,u\big\ran\\
\ns\ds=\big\lan\F_2u(\cd),u(\cd)\big\ran+2\big\lan u(\cd),\f_1\big\ran+\f_0,\ea\ee
with
$$\left\{\2n\ba{ll}
\ns\ds\F_2=\h F_1^*\cG\h F_1+F_1^*\cQ F_1+\cS F_1+F_1^*\cS^*+\cR,\\
\ns\ds\f_1=\h F_1^*\[\cG(\h F_0x+\h f_0)+g\]+F_1^*\[\cQ(F_0x+f_0)+q\]+\cS(F_0x+f_0)+\rho,\\
\ns\ds\f_0=\lan\cG(\h F_0x+\h f_0),\h F_0x+\h f_0\ran+2\lan g,\h F_0x+\h f_0\ran+\lan\cQ(F_0x+f_0),F_0x+f_0\ran+2\lan q, F_0x+f_0\ran.\ea\right.$$
Clearly,
$$\F_2\in\sS\(\cU[t,T]\),\qq\f_1\in\cU[t,T],\qq\f_0\in\dbR.$$
Therefore, according to \cite{Mou-Yong 2006}, we have the following result.

\bp{Mou-Yong} \sl For any $(t,x)\in\sD$, the map $u(\cd)\mapsto J(t,x;u(\cd))$ admits a minimum in
$\cU[t,T]$ if and only if
\bel{F_2>0}\F_2\ges0,\qq\f_1\in\sR\big(\F_2\big)\equiv\hb{ the range of $\F_2$}.\ee
In particular, if the following holds:
\bel{F_2>d}\F_2\ges\d I,\ee
for some $\d>0$, then \rf{F_2>0} holds and the map $u(\cd)\mapsto J(t,x;u(\cd))$ admits a unique minimum given by the following:
\bel{bar u}\bar u(\cd)=-\F_2^{-1}\f_1(\cd).\ee

\ep

It is clear that if
\bel{R>0}\cG\ges0,\q\cQ(\cd)-\cS(\cd)^*\cR(\cd)^{-1}\cS(\cd)\ges0,\q\cR(\cd)\ges\d I,\ee
then \rf{F_2>d} holds. Thus, under (H1)--(H2) and \rf{F_2>d}, Problem (OLQ) admits a unique open-loop optimal control.

\ms

If we denote
$$J(t,x;u(\cd))=J(t,x;u(\cd);b(\cd),\si(\cd),g,q(\cd),\rho(\cd)),$$
indicating the dependence on the nonhomogeneous terms $b(\cd),\si(\cd)$ in the state equation and the linear weighting coefficients $g,q(\cd),\rho(\cd)$ in the cost functional, then we define
\bel{J0}J^0(t;u(\cd))=J(t,0;u(\cd);0,0,0,0,0)=\lan\F_2u(\cd),u(\cd)\ran,\qq\forall u(\cd)\in\cU[t,T].\ee
Hence, we see that the following is true.

\bp{convexity equivalence} \sl Let {\rm (H1)--(H2)} hold. Then the following are equivalent:

\ms

{\rm(i)} $\F_2\ges0$;

\ms

{\rm(ii)} $J^0(t;u(\cd))\ges0$ for all $u(\cd)\in\cU[t,T]$;

\ms

{\rm(iii)} $u(\cd)\mapsto J(t,x;u(\cd))$ is convex.

\ep

\subsection{BSDEs with operator coefficients}

In this subsection, let us consider the following BSDE with operator coefficients:
\bel{Sec2_BSDE}\left\{\2n\ba{ll}
\ns\ds dY(s)=-\(\cA(s)^*Y(s)+\cC(s)^*Z(s)+\f(s)\)ds+Z(s)dW(s),\q s\in [t,T],\\
\ns\ds Y(T)=\xi\in L^2(\O;\dbR^n).\ea\right.\ee
A pair $(Y(\cd),Z(\cd))$ of adapted processes is called an adapted solution if the above is satisfied in the usual It\^o's sense. We have the following well-posedness and regularity result for the above BSDE.

\bp{BSDE well-posedness} \sl Suppose {\rm(H1)} holds and $\f(\cd)\in L_\dbF^2(\O;L^1(t,T;\dbR^n))$.  Then \rf{Sec2_BSDE} admits a unique adapted solution $(Y(\cd),Z(\cd))\in L_\dbF^2(\O;C([t,T];\dbR^n))\times L_\dbF^2(t,T;\dbR^n)$, and the following estimate holds:
\bel{Sec2_BSDE_Esup}\dbE\[\sup_{s\in[t,T]}|Y(s)|^2+\int_t^T|Z(s)|^2ds\]\les
K\[\|\xi\|_2^2+\Big\| \int_t^T |\f(s)| ds \Big\|_2^2\],\ee
where $K>0$ is a constant depending on $\|\cA(\cd)\|_1$ and $\|\cC(\cd)\|_2$.

\ep

\rm

\it Proof. \rm Denote
$$\cM[t,T]=L_\dbF^2(\O;C([t,T];\dbR^n))\times L_\dbF^2(t,T;\dbR^n),$$
and
$$\|(Y(\cd),Z(\cd))\|_{\cM[t,T]}=\Big\{\dbE\[\sup_{s\in[t,T]}|Y(s)|^2\]
+\dbE\int_t^T|Z(s)|^2ds\Big\}^{1\over2}.$$
It is clear that $\|\cd\|_{\cM[t,T]}$ is a norm under which $\cM[t,T]$ is a Banach space. We introduce a map $\sT:\cM[t,T]\to\cM[t,T]$ by the following:
$$\sT(y(\cd),z(\cd))=(Y(\cd),Z(\cd)),$$
where $(y(\cd),z(\cd))\in\cM[t,T]$ and $(Y(\cd),Z(\cd))$ is the adapted solution to the following BSDE:
$$Y(s)=\xi+\int_s^T\(\cA(r)^*y(r)+\cC(r)^*z(r)+\f(r)\)dr-\int_s^TZ(r)dW(r),\q s\in[t,T].$$
Then a standard result for BSDEs leads to the following estimate:
$$\ba{ll}
\ns\ds\dbE\[\sup_{r\in[s,T]}|Y(r)|^2+\int_s^T|Z(r)|^2dr\]\les K\Big\{\|\xi\|_2^2+\Big\|\int_s^T\(|\cA(r)^*y(r)|+|\cC(r)^*z(r)|+|\f(r)|\)dr\Big\|_2^2\Big\}\\
\ns\ds\les K\Big\{ \|\xi\|_2^2 +\Big\| \int_s^T |\cA(r)^*y(r)|dr \Big\|_2^2 +\Big\| \int_s^T |\cC(r)^*z(r)| dr \Big\|_2^2 +\Big\| \int_s^T |\f(r)|dr \Big\|_2^2 \Big\} \\
\ns\ds\les K\Big\{\|\xi\|_2^2+\( \int_s^T \|\cA(r)^* y(r)\|_2 dr \)^2 +\(\int_s^T \|\cC(r)^*z(r)\|_2dr\)^2 +\Big\| \int_s^T |\f(r)|dr \Big\|_2^2\Big\}\\
\ns\ds\les\1n K\Big\{\|\xi\|_2^2\1n+\1n\(\int_s^T\2n\|\cA(r)\|dr\)^2\dbE\[\sup_{r\in[s,T]}|y(r)|^2\]
\2n+\2n\(\int_s^T\2n\|\cC(r)\|^2dr\)\(\1n\int_s^T\2n\|z(r)\|_2^2dr\1n\)\1n+\1n\Big\| \int_s^T |\f(r)|dr \Big\|_2^2\Big\}\\
\ns\ds\les\1n K\Big\{\|\xi\|_2^2\1n+\1n\Big\| \int_s^T |\f(r)|dr
\Big\|_2^2\1n+\1n\[\(\int_s^T\2n\|\cA(r)\|dr\)^2
\2n+\2n\(\int_s^T\2n\|\cC(r)\|^2dr\)\]\[\dbE\(\sup_{r\in[s,T]}|y(r)|^2\)\1n+\1n
\int_s^T\2n\|z(r)\|_2^2dr\]\Big\}.\ea$$
The above estimate implies that the map $(y(\cd),z(\cd))\mapsto(Y(\cd),Z(\cd))$ is a contraction on $\cM[T-\d,T]$ for some small $\d>0$. Then a routine argument gives the well-posedness of BSDE \rf{Sec2_BSDE} on $[t,T]$, and estimate \rf{Sec2_BSDE_Esup} holds. \endpf

\section{Open-Loop Optimal Controls and FBSDEs}

\subsection{Solvability of Problem (OLQ)}

Let us first give the following definition.

\bde{Sec3_Def_Optim_Con} \rm A control process $\bar u(\cd)\in\cU[t,T]$ is called an {\it open-loop optimal control}
of Problem (OLQ) at $(t,x)\in\sD$ if
\bel{inf J(u)}J(t,x;\bar u(\cd))=\inf_{u(\cd)\in\cU[t,T]}J(t,x;u(\cd)).\ee
If $\bar u(\cd)\in\cU[t,T]$ exists satisfying \rf{inf J(u)}, we say
that Problem (OLQ) is {\it open-loop solvable} at $(t,x)\in\mathscr
D$, and $\bar X(\cd)\equiv X(\cd;t,x,\bar u(\cdot))$ is called the
{\it open-loop optimal state process}. \ede

The following gives a characterization of optimal open-loop control of Problem (OLQ).

\bt{optimal-condition} \sl Let {\rm(H1)--(H2)} hold. Given
$(t,x)\in\sD$, then $\bar u(\cd)\in\cU[t,T]$ is an open-loop optimal
control of Problem (OLQ) at $(t,x)\in\sD$ with $\bar X(\cd)$ being
the corresponding open-loop optimal state process if and only if
$u(\cd)\mapsto J(t,x;u(\cd))$ is convex and the following FBSDE with
operator coefficients:
\bel{FBSDE1}\left\{\2n\ba{ll}
\ns\ds\3n\ba{ll}d\bar X(s)=\(\cA(s)\bar X(s)+\cB(s)\bar u(s)+b(s)\)ds+\(\cC(s)\bar X(s)+\cD(s)\bar u(s)+\si(s)\)dW(s),\\
\ns\ds d\bar Y(s)=-\(\cA(s)^*\bar Y(s)\1n+\1n\cC(s)^*\bar Z(s)\1n+\1n\cQ(s)\bar X(s)\1n+\1n\cS(s)^*\bar u(s)\1n+\1n q(s)\)ds\1n+\1n \bar Z(s)dW(s),\ea\q s\in[t,T],\\
\ns\ds\bar X(t)=x,\qq \bar Y(T)=\cG\bar X(T)+g,\ea\right.\ee
admits an adapted solution $(\bar X(\cd),\bar Y(\cd),\bar Z(\cd))$ such that the following stationarity condition holds:
\bel{Ru+...=0}\cR(s)\bar u(s)+\cB(s)^*\bar Y(s)+\cD(s)^*\bar Z(s)+\cS(s)\bar X(s)+\rho(s)=0,\qq s\in[t,T],~\as\ee

\et

\it Proof. \rm For $(t,x)\in\sD$ and $u(\cd),\bar
u(\cd)\in\cU[t,T]$, let $X(\cd)=X(\cd\,;t,x,u(\cd))$ and $\bar
X(\cd)=X(\cd\,;t,x,\bar u(\cd))$ be the state process \rf{state}
corresponding to $u(\cd)$ and $\bar u(\cd)$, respectively. Denote
$$\h X(\cd)=X(\cd)-\bar X(\cd),\qq \h u(\cd)=u(\cd)-\bar u(\cd).$$
Then $\h X(\cd)$ satisfies the following FSDE:
$$\left\{\2n\ba{ll}
\ns\ds d\h X(s)=\(\cA(s)\h X(s)+\cB(s)\h u(s)\)ds+\(\cC(s)\h X(s)+\cD(s)\h u(s)\)dW(s),\qq s\in[t,T],\\
\ns\ds\h X(t)=0.\ea\right.$$
Now, let $(\bar Y(\cd),\bar Z(\cd))$ be the adapted solution to the BSDE in \rf{FBSDE1}. Then we have the following duality:
$$\ba{ll}
\ns\ds\dbE\lan\cG\bar X(T)+g,\h X(T)\ran=\dbE\lan\h X(T),\bar Y(T)\ran\\
\ns\ds=\dbE\int_t^T\(\lan\cA(s)\h X(s)+\cB(s)\h u(s),\bar Y(s)\ran+\lan\cC(s)\h X(s)+\cD(s)\h u(s),\bar Z(s)\ran\\
\ns\ds\qq\qq-\lan\h X(s),\cA(s)^*\bar Y(s)+\cC(s)^*\bar Z(s)+\cQ(s)\bar X(s)+\cS(s)^*\bar u(s)+q(s)\ran \)ds\\
\ns\ds=\dbE\int_t^T\(-\lan\h X(s),\cQ(s)\bar X(s)+\cS(s)^*\bar
u(s)+q(s)\ran+\lan\h u(s),\cB(s)^*\bar Y(s)+\cD(s)^*\bar Z(s)\ran\)ds.\ea$$
Hence,
$$\ba{ll}
\ns\ds J(t,x;u(\cd))-J(t,x;\bar u(\cd)) =\dbE\Big\{\lan\cG[X(T)+\bar X(T)],\h X(T)\ran+2\lan g,\h X(T)\ran\\
\ns\ds\qq+\int_t^T\(\lan\cQ(s)[X(s)+\bar X(s)],\h X(s)\ran+\lan\cS(s)[X(s)+\bar X(s)],\h u(s)\ran +\lan\cS(s)\h
X(s),u(s)+\bar u(s)\ran\\
\ns\ds\qq\qq\q+\lan\cR(s)[u(s)+\bar u(s)],\h u(s)\ran +2\lan q(s),\h X(s)\ran +2\lan\rho(s),\h u(s)\ran\)ds\Big\}\\
\ns\ds=\dbE\Big\{\lan\cG\h X(T),\h X(T)\ran+2\lan\cG\bar X(T)+g,\h X(T)\ran\\
\ns\ds\qq+\int_t^T\(\lan\cQ(s)\h X(s),\h X(s)\ran+2\lan\cS(s)\h X(s),\h u(s)\ran +\lan\cR(s)\h u(s),
\h u(s)\ran\\
\ns\ds\qq\qq\q+2\[\lan\cQ(s)\bar X(s),\h X(s)\ran+\lan\cS(s)\bar
X(s),\h u(s)\ran
+\lan\cS(s)\h X(s),\bar u(s)\ran\\
\ns\ds\qq\qq\q+\lan\cR(s)\bar u(s),\h u(s)\ran+\lan q(s),\h X(s)\ran+\lan\rho(s),\h u(s)\ran\]\)ds\Big\}\\
\ns\ds=\dbE\Big\{\lan\cG \hat X(T), \hat X(T)\ran+\int_t^T\(\lan\cQ(s) \hat X(s), \hat X(s)\ran +2\lan\cS(s)\hat X(s), \hat u(s)\ran+\lan\cR(s) \hat u(s), \hat u(s)\ran\)ds\\
\ns\ds\qq\qq\q+2\[\lan\cG\bar X(T)+g, \hat X(T)\ran+\int_t^T\(\lan\cQ(s)\bar X(s), \hat X(s)\ran +\lan\cS(s)\bar X(s), \hat u(s)\ran\\
\ns\ds\qq\qq\q+\lan\cS(s) \hat X(s),\bar u(s)\ran +\lan\cR(s)\bar u(s), \hat u(s)\ran+\lan q(s), \hat X(s)\ran+\lan\rho(s), \hat u(s)\ran\)ds\]\Big\}\\
\ns\ds=J^0(t;\hat u(\cd))+2\dbE\int_t^T \Big\lan\cR(s)\bar
u(s)+\cB(s)^*\bar Y(s)+\cD(s)^*\bar Z(s)+\cS(s)\bar X(s)+\rho(s),
\hat u(s)\Big\ran ds.\ea$$
Consequently,
$$\ba{ll}
\ns\ds J(t,x;\bar u(\cd)+\a u(\cd))-J(t,x;\bar u(\cd))=\a^2J^0(t;u(\cd)-\bar u(\cd))\\
\ns\ds\qq\qq\qq\qq+2\a\dbE \int_t^T\Big\lan\cR(s)\bar
u(s)+\cB(s)^*\bar Y(s)+\cD(s)^*\bar Z(s)+\cS(s)\bar
X(s)+\rho(s),u(s)-\bar u(s)\Big\ran ds.\ea$$
Thus, if constraint \rf{Ru+...=0} holds and $u(\cd)\mapsto J(t,x;u(\cd))$ is convex (which is equivalent to $J^0(t;u(\cd)-\bar u(\cd))\ges0$, see Proposition \ref{convexity equivalence}), we have
the optimality of $\bar u(\cd)$. Conversely, if $\bar u(\cd)\in\cU[t,T]$ is a minimum of $u(\cd)\mapsto J(t,x;u(\cd))$, then by letting $\a\to\i$ in the above, we see that $u(\cd)\mapsto J(t,x;u(\cd))$ is convex and the constraint \rf{Ru+...=0} holds. \endpf

\ms

We note that \rf{FBSDE1}--\rf{Ru+...=0} is a coupled linear FBSDE
with operator coefficients. The above theorem tells us that the
open-loop solvability of Problem (OLQ) is equivalent to the
solvability of an FBSDE with operator coefficients. A similar result
for LQ problem with constant (matrix) coefficients were established in
\cite{Sun-Yong 2016}. The proof presented above is similar to that
found in \cite{Sun-Yong 2016}, with some simplified arguments.

\ms

\subsection{Well-posedness of FBSDEs with operator coefficients}

We now look at the solvability of FBSDE \rf{FBSDE1}--\rf{Ru+...=0}. To abbreviate the notations, we drop the bars in $\bar{X},\bar{Y},\bar{Z},\bar{u}$ of \rf{FBSDE1} and \rf{Ru+...=0},
that is, we consider the following:
\bel{Sec5_FBSDE}\left\{\2n\ba{ll}
\ns\ds dX(s)=\(\cA(s)X(s)+\cB(s)u(s)+b(s)\)ds+\(\cC(s)X(s)+\cD(s)u(s)+\si(s)\)dW(s),\\
\ns\ds dY(s)=-\(\cA(s)^*Y(s)+\cC(s)^*Z(s)+\cQ(s)X(s)+\cS(s)^*u(s)+q(s)\)ds+Z(s)dW(s),\\
\ns\ds X(t)=x,\qq Y(T)=\cG X(T)+g,\\
\ns\ds\cR(s)u(s)+\cB(s)^*Y(s)+\cD(s)^*Z(s)+\cS(s)X(s)+\rho(s)=0.\ea\right.s\in [t,T],\ee
For the well-posedness of the above equation, we need the following hypothesis.

\ms

{\bf(H3)} 
%
%
For some $\d>0$,
$$\cG\ges 0,\q\cQ(s)-\cS(s)^*\cR(s)^{-1}\cS(s)\ges0,\q\cR(s)\ges\d I,\qq s\in[0,T],$$
where $I$ in the above is the identity operator on $L^2(\O;\dbR^m)$.

\ms

Note that the last condition in (H3) ensures the existence of the inverse $\cR^{-1}$ of operator
$\cR$, and
\bel{R-1}\dbE\lan \cR(s)^{-1}u,u\ran\ges{\d\over\|\cR(s)\|^2}
\dbE|u|^2, \hb{ for any } u\in L^2_{\cF_s} (\O;\dbR^m),\ \ae, \ s\in
[t,T].\ee
On the other hand, by Schur's Lemma, one sees that (H3) is
equivalent to the following:
$$\cG\ges0,\qq\begin{pmatrix}\cQ(\cd) & \cS(\cd)^*\\
                             \cS(\cd) & \cR(\cd)\end{pmatrix}\ges 0,\q\cR(\cd)\ges\d I.$$

\ms

Now, under (H3), $u(\cd)$ can be expressed explicitly as
\bel{Sex5_u_explicit}
u(s)=-\cR(s)^{-1}\(\cB(s)^*Y(s)+\cD(s)^*Z(s)+\cS(s)X(s)+\rho(s)\),\q s\in [t,T].
\ee
By substituting it into the first two equations of \rf{Sec5_FBSDE}, the system becomes
a coupled FBSDEs with operator coefficients as follows (suppressing $s$):
\bel{sec5_coupled_case}\left\{\2n\ba{ll}
\ns\ds dX=\((\cA -\cB\cR^{-1}\cS) X-\cB\cR^{-1}\cB^*Y-\cB\cR^{-1}\cD^*Z-\cB\cR^{-1} \rho+b\)ds\\
\ns\ds \qq\q +\((\cC-\cD\cR^{-1}\cS) X-\cD\cR^{-1}\cB^*Y-\cD\cR^{-1}\cD^*Z-\cD\cR^{-1}\rho+\si\)dW(s),\\
\ns\ds dY=-\((\cQ -\cS^*\cR^{-1}\cS)X+(\cA^*-\cS^*\cR^{-1}\cB^*)Y+(\cC^*-\cS^*\cR^{-1} \cD^*)Z\\
\ns\ds \qq\qq -\cS^*\cR^{-1} \rho+q\)ds +ZdW(s),\\
\ns\ds X(t)=x,\qq Y(T)=\cG X(T)+g.\ea\right.\ s\in [t,T],\ee

\br{}\rm From {\rm(H1), (H2)}, it is easy to check that
$$\ba{ll}
\ns\ds\cA(\cd)-\cB(\cd)\cR(\cd)^{-1}\cS(\cd)\in\sL^1_\dbF(L^2(\O;\dbR^n)),\\
\ns\ds\cC(\cd)-\cD(\cd)\cR(\cd)^{-1}\cS(\cd)\in\sL^2_\dbF(L^2(\O;\dbR^n)),\\
\ns\ds\cQ(\cd)-\cS(\cd)^*\cR(\cd)^{-1}\cS(\cd)\in\sS^1_\dbF(L^2(\O;\dbR^n)).
\ea$$
For example,
$$\ba{ll}
\ns\ds \int_0^T\|\cA(s)-\cB(s)\cR(s)^{-1}\cS(s)\|ds\les\int_0^T\(\|\cA(s)\|+\|\cB(s)\cR(s)^{-1}\cS(s)\|\)ds\\
\ns\ds \les\int_0^T\|\cA(s)\|ds+\(\int_0^T\|\cB(s)\|^2ds\)^{1\over2}\sup_{s\in[0,T]}\|\cR(s)^{-1}\|\(\int_0^T\|\cS(s)\|^2ds\)^{1\over2}<\i.
\ea$$
The remaining two can be verified similarly.

\er

The following is the main result of this subsection.

\bt{Sec5_wellposedness} \sl Let {\rm(H1)--(H3)} hold. Then there
exists a unique adapted solution $(X(\cd),Y(\cd),Z(\cd))\in
[L^2_\dbF(\O;C([t,T];\dbR^n))]^2 \times L^2_\dbF(t,T;\dbR^n)$ to the
coupled FBSDE \rf{sec5_coupled_case} with operator coefficients, and
the following estimate holds:
\bel{Yu_FBSDE_Estimate}\2n\ba{ll}
\ns\ds\dbE\[\sup_{s\in [t,T]}|X(s)|^2+\sup_{s\in [t,T]} |Y(s)|^2+ \int_t^T|Z(s)|^2ds\]\\
\ns\ds\les K\Big\{\|x\|_2^2+\|g\|_2^2+\Big\|\int_t^T |b(s)|
ds\Big\|_2^2 +\Big\| \int_t^T |q(s)|ds
\Big\|_2^2+\int_t^T\|\sigma(s)\|_2^2ds +\int_t^T \|\rho(s)\|_2^2 ds\Big\}
,\ea\ee
where $K>0$ is a constant depending on $\|\cA(\cd)\|_1$,
$\|\cB(\cd)\|_2$, $\|\cC(\cd)\|_2$, $\|\cD(\cd)\|_\i$,
$\|\cQ(\cd)\|_1$, $\|\cS(\cd)\|_2$, $\|\cR(\cd)^{-1}\|_\i$ and
$\|\cG\|$. \et

Before presenting a proof of Theorem \ref{Sec5_wellposedness}, we
introduce the following auxiliary FBSDE with operator coefficients
parameterized by $\a\in [0,1]$:
\bel{Sec5_FBSDE_alpha}\3n\3n\3n\2n_\a\qq\left\{
\2n\ba{ll}
\ns\ds dX^\a=\(\a(\cA-\cB\cR^{-1}\cS)X^\a-\cB\cR^{-1}\cB^*Y^\a-\cB\cR^{-1}\cD^*Z^\a+\f\)ds\\
\ns\ds\qq\qq+\(\a(\cC-\cD\cR^{-1}\cS)X^\a-\cD\cR^{-1}\cB^*Y^\a-\cD\cR^{-1}\cD^*Z^\a)+\psi\)dW(s),\\
\ns\ds dY^{\a}=-\(\a(\cQ -\cS^*\cR^{-1}\cS)X^\a+\a(\cA^*-\cS^*\cR^{-1}\cB^*)Y^\a\\
\ns\ds\qq\qq+\a(\cC^*-\cS^*\cR^{-1} \cD^*)Z^\a+\g\)ds +Z^\a dW(s), \\
\ns\ds X^\a(t)=\xi,\qq Y^\a(T)=\a\cG X^\a(T)+\eta,\ea\qq s\in[t,T],\right.\ee
where
$$(\xi,\f,\psi,\g,\eta)\in\dbM[t,T]\equiv L^2_{\cF_t}(\O;\dbR^n)\times
L^2_\dbF(\O;L^1(t,T;\dbR^n))\times L^2_\dbF(t,T;\dbR^n)\times
L^2_\dbF(\O;L^1(t,T;\dbR^n))\times L^2(\O;\dbR^n).$$
The following gives an {\it a priori} estimate of the adapted solution $(X^\a(\cd),Y^\a(\cd),Z^\a(\cd))$ to the parameterized equation \rf{Sec5_FBSDE_alpha}$_\a$.

\bp{Sec5_priori_estimate}
\sl Let {\rm(H1)--(H3)} hold and $\a\in[0,1]$ be given. Let
$(X_1,Y_1,Z_1)$ and $(X_2,Y_2,Z_2)$ be the solutions to
\rf{Sec5_FBSDE_alpha}$_\a$ corresponding to the different
$(\xi_1,\f_1,\psi_1,\g_1,\eta_1)$,
$(\xi_2,\f_2,\psi_2,\g_2,\eta_2)\in\dbM[t,T]$, respectively. Then
the following estimate holds:
\bel{Sec5_ee0}\2n\ba{ll}
\ns\ds\dbE\[\sup_{s\in [t,T]}|X_1(s)-X_2(s)|^2+\sup_{s\in [t,T]} |Y_1(s)-Y_2(s)|^2+ \int_t^T|Z_1(s)-Z_2(s)|^2ds\]\\
\ns\ds\les K\Big\{\|\xi_1-\xi_2\|_2^2+\|\eta_1-\eta_2\|_2^2+\Big\|\int_t^T |\f_1(s)-\f_2(s)| ds\Big\|_2^2+\int_t^T\|\psi_1(s)-\psi_2(s)\|_2^2ds \\
\ns\ds\qq\qq +\Big\| \int_t^T |\g_1(s)-\g_2(s)|ds
\Big\|_2^2\Big\},\ea\ee
where $K>0$ is a constant depending on $\|\cA(\cd)-\cB(\cd)\cR(\cd)^{-1}\cS(\cd)\|_1$, $\|\cC(\cd)-\cD(\cd)\cR(\cd)^{-1}\cS(\cd)\|_2$,
$\|\cQ(\cd)-\cS(\cd)^*\cR(\cd)^{-1}\cS(\cd)\|_1$, $\|\cB(\cd)\|_2$, $\|\cD(\cd)\|_\i$, $\|\cR(\cd)^{-1}\|_\i$, $\|\cG\|$, and independent of $\a\in[0,1]$.

\ep

\rm

\it Proof. \rm For convenience, we denote $(\h X,\h Y,\h Z):=(X_1-X_2,Y_1-Y_2,Z_1-Z_2)$, which
satisfies
\bel{Sec5_FBSDE_difference}\left\{\2n\ba{ll}
\ns\ds d\h X=\(\a(\cA-\cB\cR^{-1}\cS)\h X-\cB\cR^{-1}\cB^*\h Y-\cB\cR^{-1}\cD^*\h Z+\h\f\)ds\\
\ns\ds \qq\q+\(\a(\cC-\cD\cR^{-1}\cS)\h X-\cD\cR^{-1}\cB^*\h Y-\cD\cR^{-1}\cD^*\h Z+\h\psi\)dW(s), \q s\in [t,T], \\
\ns\ds d\h Y=-\(\a(\cQ -\cS^*\cR^{-1}\cS)\h X+\a(\cA^*-\cS^*\cR^{-1}\cB^*)\h Y\\
\ns\ds\qq\qq+\a(\cC^*-\cS^*\cR^{-1}\cD^*)\h Z+\h\g\)ds+\h Z dW(s), \q s\in [t,T],\\
\ns\ds\h X(t)=\h\xi,\qq \h Y(T)=\a\cG\h X(T)+\h\eta,\ea\right.\ee
with
$(\h\xi,\h\f,\h\psi,\h\g,\h\eta)=(\xi_1-\xi_2,\f_1-\f_2,\psi_1-\psi_2,\g_1-\g_2,\eta_1-\eta_2)$.
\ms

For the forward equation, by applying Proposition
\ref{Well-posedness-state}, we have
$$\ba{ll}
\ns\ds\dbE\Big[\sup_{s\in [t,T]}|\h X(s)|^2\Big]\\
\ns\ds\les K\Big\{\|\h\xi\|_2^2+\Big\|\int_t^T |-\cB\cR^{-1}(\cB^*\h
Y+\cD^*\h Z)
+\h\f | dr\Big\|_2^2+\int_t^T\|-\cD\cR^{-1}(\cB^*\h Y +\cD^*\h Z)+\h\psi\|_2^2dr\Big\}\\
\ns\ds\les K\Big\{\|\h\xi\|_2^2+\Big\| \int_t^T |\h \f| dr \Big\|_2^2+\int_t^T\|\h\psi\|_2^2dr\\
\ns\ds\qq\q+\(\|\cB(\cd)\|_2^2+\|\cD(\cd)\|_\i^2\)\|\cR(\cd)^{-1}\|_\i^2\int_t^T\|\cB^*\h Y+\cD^*\h Z\|_2^2dr\Big\}\\
\ns\ds\les K\Big\{\|\h\xi\|_2^2+\Big\| \int_t^T |\h \f| dr
\Big\|_2^2+\int_t^T\|\h\psi\|_2^2dr+\int_t^T\|\cB^*\h Y +\cD^*\h
Z\|_2^2dr\Big\},\ea$$
where $K$ depends on $\|\cA(\cd)-\cB(\cd)\cR(\cd)^{-1}\cS(\cd)\|_1$, $\|\cC(\cd)-\cD(\cd)\cR(\cd)^{-1}\cS(\cd)\|_2$,
 $\|\cB(\cd)\|_2$, $\|\cD(\cd)\|_\i$, $\|\cR(\cd)^{-1}\|_\i$, independent of $\a$.

\ms

Similarly, for the backward equation, Proposition \ref{BSDE well-posedness} leads to
$$\ba{ll}
\ns\ds\dbE\[\sup_{s\in[t,T]}|\h Y(s)|^2+\int_t^T|\h Z(s)|^2ds\]\\
\ns\ds\les K\Big\{\|\a\cG\h X(T)+\h\eta\|_2^2+\Big\|\int_t^T |\a(\cQ -\cS^*\cR^{-1}\cS)\h X +\hat\g | ds\Big\|_2^2\bigg\}\\
\ns\ds\les K\Big\{\|\h\eta\|_2^2+\Big\| \int_t^T |\h\g| ds \Big\|_2^2+\(\|\cG\|^2+\|\cQ(\cd)-\cS(\cd)^*\cR(\cd)^{-1}\cS(\cd)\|_1^2\)\dbE[\sup_{s\in [t,T]} |\h X(s)|^2\]\Big\}\\
\ns\ds\les K\Big\{\|\h\eta\|_2^2+\Big\| \int_t^T |\h\g| ds
\Big\|_2^2+\dbE\[\sup_{s\in [t,T]}|\h X(s)|^2\]\Big\},\ea$$
where $K$ depends on $\|\cA(\cd)-\cB(\cd)\cR(\cd)^{-1}\cS(\cd)\|_1$, $\|\cC(\cd)-\cD(\cd)\cR(\cd)^{-1}\cS(\cd)\|_2$,
$\|\cQ(\cd)-\cS(\cd)^*\cR(\cd)^{-1}\cS(\cd)\|_1$, $\|\cG\|$, independent of $\a$. Then, we get
\bel{Sec3.2_Eq1}\ba{ll}
\ns\ds\dbE\[\sup_{s\in [t,T]}|\h X(s)|^2+\sup_{s\in [t,T]}|\h Y(s)|^2+\int_t^T|\h Z(s)|^2ds\]\\
\ns\ds\les K\Big\{\|\h\xi\|_2^2+\|\h\eta\|_2^2+\Big\| \int_t^T
|\h\f| dr \Big\|_2^2+\int_t^T\|\h\psi\|_2^2dr+\Big\| \int_t^T |\h\g|
dr \Big\|_2^2+\int_t^T\|\cB^*\h Y+\cD^*\h Z\|_2^2dr\Big\}.\ea\ee
with $K$ depending on $\|\cA(\cd)-\cB(\cd)\cR(\cd)^{-1}\cS(\cd)\|_1$, $\|\cC(\cd)-\cD(\cd)\cR(\cd)^{-1}\cS(\cd)\|_2$,
$\|\cQ(\cd)-\cS(\cd)^*\cR(\cd)^{-1}\cS(\cd)\|_1$, $\|\cB(\cd)\|_2$, $\|\cD(\cd)\|_\i$, $\|\cR(\cd)^{-1}\|_\i$, $\|\cG\|$, independent of $\a$.

\ms

On the other hand, by applying It\^o's formula to $\lan
\h X(\cd),\h Y(\cd)\ran $, we have
$$\ba{ll}
\ns\ds\dbE\Big\{\lan\h X(T),\a\cG\h X(T)+\h\eta\ran-\lan\h\xi,\h Y(t)\ran\Big\}=\dbE\Big\{\lan\h X(T),\h Y(T)\ran-\lan\h X(t),\h Y(t)\ran\Big\}\\
\ns\ds=\dbE\int_0^T\Big\{\lan\a(\cA-\cB\cR^{-1}\cS)\h X-\cB\cR^{-1}\cB^*\h Y-\cB\cR^{-1}\cD^*\h Z+\h\f,\h Y\ran\\
\ns\ds\qq\qq-\lan\h X,\a(\cQ-\cS^*\cR^{-1}\cS)\h X+\a(\cA^*-\cS^*\cR^{-1}\cB^*)\h Y+\a(\cC^*-\cS^*\cR^{-1}\cD^*)\h Z+\h\g\ran\\
\ns\ds\qq\qq+\lan\a(\cC-\cD\cR^{-1}\cS)\h X-\cD\cR^{-1}\cB^*\h Y-\cD\cR^{-1}\cD^*\h Z+\h\psi,\h Z\ran\Big\}ds\\
\ns\ds=\dbE\int_0^T\Big\{\lan-\cB\cR^{-1}\cB^*\h Y-\cB\cR^{-1}\cD^*\h Z,\h Y\ran-\a\lan\h X,(\cQ-\cS^*\cR^{-1}\cS)\h X\ran\\
\ns\ds\qq\qq+\lan-\cD\cR^{-1}\cB^*\h Y-\cD\cR^{-1}\cD^*\h Z,\h Z\ran+\lan\h\f,\h Y\ran+\lan\h\psi,\h Z\ran-\lan\h X,\h\g\ran\Big\}ds\\
\ns\ds=-\dbE\int_0^T\Big\{\lan\cR^{-1}(\cB^*\h Y+\cD^*\h Z),\cB^*\h Y+\cD^*\h Z\ran-\a\lan\h X,(\cQ-\cS^*\cR^{-1}\cS)\h X\ran+\lan\h\f,\h Y\ran+\lan\h\psi,\h Z\ran-\lan\h X,\h\g\ran\Big\}ds.\ea$$
Hence,
\bel{Sec3_monotonicity}\ba{ll}
\ns\ds\dbE\int_t^T\lan\cR^{-1}(\cB^*\h Y+\cD^*\h Z),\cB^*\h Y+\cD^*\h Z\ran dr\\
\ns\ds=-\a\dbE\Big\{\lan\cG\h X(T),\h X(T)\ran+\int_t^T\lan(\cQ-\cS^*\cR^{-1}\cS)\h X,\h X\ran dr\Big\}\\
\ns\ds\q+\dbE\Big\{\lan\h\xi,\h Y(t)\ran-\lan\h X(T),\h\eta\ran+\int_t^T\[\lan\h\f,
\h Y\ran+\lan\h\psi,\h Z\ran-\lan\h\g,\h X\ran\]dr\Big\}.\ea\ee
By (H3), we know
$$\dbE\int_t^T|\cB^*\h Y+\cD^*\h Z|^2ds\les\frac{\|\cR\|^2}{\d}\dbE\Big\{\lan\h\xi,\h Y(t)\ran
-\lan\h\eta,\h X(T)\ran+\int_t^T\[\lan\h\f,\h Y\ran+\lan\h\psi,\h Z\ran
-\lan\h\g,\h X\ran\]dr\Big\}.$$
Then, for any $\e> 0$, we have
\bel{Sec3.2_Eq2}\ba{ll}
\ns\ds\dbE\int_t^T|\cB^*\h Y+\cD^*\h Z|^2ds\les\e\dbE\[\sup_{s\in[t,T]}|\h X(s)|^2+\sup_{s\in [t,T]}|\h Y(s)|^2+\int_t^T|\h Z(s)|^2ds\]\\
\ns\ds\qq\qq\qq\qq\qq\q+\frac{\|\cR\|^4}{4\e\d^2}\Big\{\|\h\xi\|_2^2+\|\h\eta\|_2^2+\Big\|
\int_t^T |\h\f| dr \Big\|_2^2 +\Big\| \int_t^T |\h\g| dr
\Big\|_2^2+\int_t^T\|\h\psi\|_2^2dr\Big\}.\ea\ee
By selecting $\e = 1/(2K)$ ($K$ is the one in \rf{Sec3.2_Eq1}),  we get the desired result
\rf{Sec5_ee0} from \rf{Sec3.2_Eq1} and \rf{Sec3.2_Eq2}.

\endpf

\br{}\rm It is known that, for general coupled FBSDEs, monotonicity
condition leads to the solvability (\cite{Hu-Peng 1995,Peng-Wu
1999}). Our condition (H3) is basically the monotonicity condition
for linear FBSDEs with operator coefficients.

\er

Next lemma gives the method of continuation.

\bl{Sec5_Contin_Lemma} \sl Let {\rm (H1)--(H3)} hold. Then there exists a constant $\e_0>0$ such that for any given $\a_0\in [0,1)$, if FBSDE \rf{Sec5_FBSDE_alpha}$_{\a_0}$ admits a unique adapted solution for any $(\xi,\f,\psi,\g,\eta)\in \dbM[t,T]$, then for $\a=\a_0+\e$ with $\e\in(0,\e_0]$, $\a_0+\e\les 1$, \rf{Sec5_FBSDE_alpha}$_\a$ also admits a unique adapted solution for any $(\xi,\f,\psi,\g,\eta)\in \dbM[t,T]$.

\el

\rm

\it Proof. \rm Let $\e_0>0$ be undetermined, and $\e\in(0,\e_0]$. We focus on the following FBSDE:
\bel{Sec5_FBSDE_in_contin}\left\{
\2n\ba{ll}
\ns\ds dX=\(\a_0(\cA-\cB\cR^{-1}\cS)X-\cB\cR^{-1}(\cB^*Y+\cD^*Z)+\e(\cA-\cB\cR^{-1}\cS)\cX+\f\)ds\\
\ns\ds\qq\q+\(\a_0(\cC-\cD\cR^{-1}\cS)X-\cD\cR^{-1}(\cB^*Y+\cD^*Z)+\e(\cC-\cD\cR^{-1}\cS)\cX +\psi\)dW(s),\\
\ns\ds dY=-\(\a_0(\cA^*-\cS^*\cR^{-1}\cB^*)Y+\a_0(\cC^*-\cS^*\cR^{-1}\cD^*)Z+\a_0(\cQ -\cS^*\cR^{-1}\cS)X\\
\ns\ds\qq\q+\e(\cA^*-\cS^*\cR^{-1}\cB^*)\cY+\e(\cC^*-\cS^*\cR^{-1}\cD^*)\cZ+\e(\cQ -\cS^*\cR^{-1}\cS)\cX+\g\)ds\\
\ns\ds\qq\q+ZdW(s), \\
\ns\ds X(t)=\xi,\qq Y(T)=\a_0\cG X(T)+\e\cG\cX(T)+\eta,\ea\right.~s\in[t,T],\ee
where  $(\cX,\cY,\cZ)\in [L_\dbF^2(\O;C([t,T];\dbR^n))]^2
\times L^2_\dbF(t,T;\dbR^n)$ is arbitrarily chosen.

\ms

Our assumption ensures the solvability of the above equation. Therefore, we can define a mapping $\cL_{\a_0+\e}$ from the space $[L_\dbF^2(\O;C([t,T];\dbR^n))]^2
\times L^2_\dbF(t,T;\dbR^n)$ into itself as follows:
$$(X(\cd),Y(\cd),Z(\cd))=\cL_{\a_0+\e} (\cX(\cd),\cY(\cd),\cZ(\cd)).$$
For another given $(\bar \cX(\cd),\bar
\cY(\cd),\bar\cZ(\cd))\in [L_\dbF^2(\O;C([t,T];\dbR^n))]^2
\times L^2_\dbF(t,T;\dbR^n)$, let
$$(\bar X(\cd),\bar Y(\cd),\bar Z(\cd))=\cL_{\a_0+\e}(\bar \cX(\cd),\bar \cY(\cd),\bar\cZ(\cd)). $$
Denote $(\h X(\cd),\h Y(\cd),\h Z(\cd))=(X(\cd)-\bar X(\cd),Y(\cd)-\bar Y(\cd),Z(\cd)-\bar
Z(\cd))$. From Proposition \ref{Sec5_priori_estimate}, we have
$$\ba{ll}
\ns\ds\dbE\[\sup_{s\in[t,T]}|\h X(s)|^2+\sup_{s\in[t,T]}|\h Y(s)|^2+\int_t^T|\h Z(s)|^2 ds\]\\
\ns\ds\les K|\e|^2\Big\{\|\cG\h X(T)\|_2^2+\(\int_t^T\|(\cA-\cB\cR^{-1}\cS)\h\cX\|_2 dr\)^2+\int_t^T\|\cC-\cD\cR^{-1}\cS)\h\cX\|_2^2dr\\
\ns\ds\qq\qq+\(\int_t^T\|(\cA^*-\cS^*\cR^{-1}\cB^*)\h\cY+(\cC^*-\cS^*\cR^{-1}\cD^*)\h\cZ +(\cQ-\cS^*\cR^{-1}\cS)\h\cX\|_2dr\)^2\Big\}\\
\ns\ds\les K|\e|^2\Big\{\|\cG\|^2\|\h X(T)\|_2^2+\(\|\cA(\cd)-\cB(\cd)\cR(\cd)^{-1}\cS(\cd)\|_1^2+\|\cC(\cd)-\cD(\cd)\cR(\cd)^{-1}\cS(\cd)
\|_2^2\\
\ns\ds\qq\qq\qq\qq\qq\qq\q +\|\cQ(\cd)-\cS(\cd)^*\cR(\cd)^{-1}\cS(\cd)\|_1^2\)\ \dbE\[\sup_{s\in [t,T]} |\h\cX(s)|^2\]\\
\ns\ds\qq\qq+\|\cA(\cd)^*-\cS(\cd)^*\cR(\cd)^{-1}\cB(\cd)^*\|_1^2\ \dbE\[\sup_{s\in
[t,T]}|\h\cY(s)|^2\]\\
\ns\ds\qq\qq+\|\cC(\cd)^*-\cS(\cd)^*\cR(\cd)^{-1}\cD(\cd)^*\|_2^2\ \dbE\int_t^T|\h\cZ(r)|^2dr\Big\}\\
\ns\ds\les K|\e|^2\dbE\[\sup_{s\in[t,T]}|\h\cX(s)|^2+\sup_{s\in [t,T]}|\h\cY(s)|^2+\int_t^T|\h
\cZ(s)|^2ds\],\ea$$
where $K>0$ is a constant independent of $\a_0$ and $\e$.
Therefore, we can choose $\e_0>0$, such that $K\e_0^2\les 1/4$.
Then, for any $\e\in[0,\e_0]$, $\cL_{\a_0+\e}$ is a
contraction map. Consequently, there exists a unique fixed point for
the mapping $\cL_{\a_0+\e}$ which is just the unique
solution to FBSDE \rf{Sec5_FBSDE_alpha}$_\a$ with $\a
=\a_0+\e$. \endpf

\ms

Now we present a proof of Theorem \ref{Sec5_wellposedness}.

\ms

\rm

\it Proof of Theorem \ref{Sec5_wellposedness}. \rm When $\a=0$,
\rf{Sec5_FBSDE_alpha}$_0$ becomes
$$\left\{\ba{ll}
\ns\ds dX^0 = \Big( -\cB\cR^{-1}(\cB^* Y^0 +\cD^* Z^0) +\f \Big) ds+\Big( -\cD\cR^{-1}(\cB^* Y^0+\cD^* Z^0) +\psi \Big) dW(s), \\
\ns\ds dY^0 = -\g ds +Z^0 dW(s),  \\
\ns\ds X^0(t) = \xi,\qq Y^0(T) =\eta,\ea\q s\in [t,T],\right.$$
whose solvability is clear. In fact, one can first solve the BSDE to obtain $(Y^0,Z^0)$, then solve the FSDE to get $X^0$.

\ms

Next, by Lemma \ref{Sec5_Contin_Lemma}, for any
$(\xi,\f,\psi,\g,\eta)\in \dbM[t,T]$, and any $\a\in[0,1]$, \rf{Sec5_FBSDE_alpha}$_\a$
is uniquely solvable. In particular, when $\a=1$, \rf{Sec5_FBSDE_alpha}$_1$ with
\begin{equation}\label{Yu_Temp1}
\xi =x,\q \f =-\cB\cR^{-1}\rho +b,\q \psi = -\cD\cR^{-1}\rho +\si,\q
\g = -\cS\cR^{-1}\rho +q,\q \eta =g
\end{equation}
becomes \rf{Sec5_FBSDE} which is also uniquely solvable.

\ms

Let $\alpha=1$, $(\xi,\varphi,\psi,\gamma,\eta)$ is given by
\eqref{Yu_Temp1}, $(\xi_2,\varphi_2,\psi_2,\gamma_2,\eta_2) =
(0,0,0,0,0)$ (the corresponding solution to FBSDE is $(0,0,0)$). By
the a priori estimate \eqref{Sec5_ee0}, we have
$$
\begin{aligned}
& \mathbb E\[\sup_{s\in [t,T]}|X(s)|^2 +\sup_{s\in [t,T]}|Y(s)|^2
+\int_t^T |Z(s)|^2 ds\]\\
& \les K\Big\{ \|x\|_2^2 +\|g\|_2^2 +\Big\| \int_t^T
|-\cB(s)\cR(s)^{-1}\rho(s) +b(s)|ds \Big\|_2^2\\
& \qquad +\int_t^T \| -\cD(s)\cR(s)^{-1}\rho(s) +\sigma(s) \|_2^2 ds
+\Big\| \int_t^T |-\cS(s)\cR(s)^{-1}\rho(s) +q(s)|ds \Big\|_2^2
\Big\}\\
& \les K\Big\{ \|x\|_2^2 +\|g\|_2^2 +\Big\|\int_t^T
|b(s)|ds\Big\|_2^2 +\int_t^T \|\sigma(s)\|_2^2 ds +\Big\|\int_t^T
|q(s)| ds\Big\|_2^2\\
& \qquad +\int_t^T \[ \|\cB(s)\cR(s)^{-1}\rho(s)\|_2^2
+\|\cD(s)\cR(s)^{-1}\rho(s)\|_2^2 +\|\cS(s)\cR(s)^{-1}\rho(s)\|_2^2
\] ds\Big\}\\
& \les K\Big\{ \|x\|_2^2 +\|g\|_2^2 +\Big\|\int_t^T
|b(s)|ds\Big\|_2^2 +\int_t^T \|\sigma(s)\|_2^2 ds +\Big\|\int_t^T
|q(s)| ds\Big\|_2^2\\
& \qquad +\(\|\cB(\cd)\|_2^2 +\|\cD(\cd)\|_\i^2 +\|\cS(\cd)\|_2^2\)
\|\cR(\cd)^{-1}\|_\i^2 \int_t^T \|\rho(s)\|_2^2 ds\Big\}.
\end{aligned}
$$
We obtain \eqref{Yu_FBSDE_Estimate}.
\endpf

\ms

\br{} \rm Note that, in our setting, the norms
$\|\cA(\cd)-\cB(\cd)\cR(\cd)^{-1}\cS(\cd)\|$,
$\|\cQ(\cd)-\cS(\cd)^*\cR(\cd)^{-1}\cS(\cd)\|$, and
$\|\cB(\cdot)\cR(\cdot)^{-1}\cB(\cdot)^*\|$ of coefficients are only
required to belong to $L^1(t,T;\dbR)$, and
$\|\cC(\cd)-\cD(\cd)\cR(\cd)^{-1}\cS(\cd)\|$,
$\|\cD(\cdot)\cR(\cdot)^{-1}\cB(\cdot)^*\|$ $\in L^2(t,T;\dbR)$.
Therefore, these five norms are not necessarily bounded. \er

\bc{Yu_Sec3.2_Corollary} \sl Under {\rm (H1)--(H3)}, Problem (OLQ) admits a unique open-loop optimal control given by \rf{Sex5_u_explicit}, where $(X(\cd),Y(\cd),Z(\cd))$ is the unique solution to FBSDE
\rf{sec5_coupled_case}.
\ec

\it Proof. \rm Let us denote by $J^0(t;u(\cd))$ the cost functional \rf{cost} when $x, b(\cd), \si(\cd), g,
q(\cd), \rho(\cd)$ are all $0$. Obviously, Assumption (H3)--(ii) implies $J^0(t;u(\cd)) \ges 0$ for all $u(\cd)\in\cU[t,T]$. By Proposition \ref{convexity equivalence}, we know $u(\cd) \mapsto J(t,x;u(\cd))$ is convex. Moreover, by Theorem \ref{Sec5_wellposedness} and the expression \rf{Sex5_u_explicit},
there exists a unique $(X(\cd), Y(\cd), Z(\cd), u(\cd))$ satisfying \rf{FBSDE1} and \rf{Ru+...=0}. Thanks to Theorem \ref{optimal-condition}, we obtain the result. \endpf

\section{Mean-Field LQ Control Problem}

We have mentioned that a major motivation of this work is the study of stochastic LQ problem of mean-field FSDE with cost functional also involving mean-field terms. In this section, we will
carry out some details for the mean-field case.

\ms

We shall use the mean-field setting \rf{sumA} and \rf{Sec2.1_Cost_Operators} given in Subsection 2.1£¬ which we rewrite here for convenience.
\bel{Sec4_Oper_in_SDE}\left\{\2n\2n\ba{ll}
\ns\ds\ba{ll}
\ns\ds\cA(s)\xi=A(s)\xi+\bar\BA(s)^\top\dbE\big[\wt\BA(s)\xi\big],\\
\ns\ds\cC(s)\xi=C(s)\xi+\bar\BC(s)^\top\dbE\big[\wt\BC(s)\xi\big],\ea\qq\forall\xi\in L^2_{\cF_s}(\O;\dbR^n),\\
\ns\ds\ba{ll}
\ns\ds\cB(s)\eta=B(s)\eta+\bar\BB(s)^\top\dbE\big[\wt\BB(s)\eta\big],\\
\ns\ds\cD(s)\eta=D(s)\eta+\bar\BD(s)^\top\dbE\big[\wt\BD(s)\eta\big],\ea\qq\forall\eta\in L^2_{\cF_s}(\O;\dbR^m),\ea\right.\ee
and
\bel{Sec4_Oper_in_cost}\left\{\2n\ba{ll}
\ns\ds\cG\xi=G\xi+\wt\BG^\top\dbE[\bar\BG\xi]+\bar\BG^\top\dbE[\wt\BG\xi]
+\wt\BG^\top\dbE[\h\BG]\dbE[\wt\BG\xi],\q\xi\in L^2(\O;\dbR^n),\\
\ns\ds\cQ(s)\xi=Q(s)\xi+\wt\BQ(s)^\top\dbE[\bar\BQ(s)\xi]+\bar\BQ(s)^\top\dbE[\wt\BQ(s)\xi]
+\wt\BQ(s)^\top\dbE[\h\BQ(s)]\dbE[\wt\BQ(s)\xi],\q\xi\in
L^2_{\cF_s}(\O;\dbR^n),\\
\ns\ds\cS(s)\xi=S(s)\xi+\wt\BR(s)^\top\dbE[\bar\BS(s)\xi]+\wt\BS(s)^\top\dbE[\wt\BQ(s)\xi]
+\wt\BR(s)^\top\dbE[\h\BS(s)]\dbE[\wt\BQ(s)\xi],\q\xi\in
L^2_{\cF_s}(\O;\dbR^n),\\
\ns\ds\cR(s)\eta=R(s)\eta+\wt\BR(s)^\top\dbE[\bar\BR(s)\eta]+\bar\BR(s)^\top\dbE[\wt\BR(s)\eta]
+\wt\BR(s)^\top\dbE[\h\BR(s)]\dbE[\wt\BR(s)\eta],\q\eta\in
L^2_{\cF_s}(\O;\dbR^m),\\
\ns\ds g=g_0+\wt\BG^\top\dbE[\bar\Bg],\q
q(s)=q_0(s)+\wt\BQ(s)^\top\dbE[\bar\Bq(s)],\q
\rho(s)=\rho_0(s)+\wt\BR(s)^\top\dbE[\bar\Brho(s)],\ea\right.\ee
with the above coefficients given by \rf{BA}--\rf{g}.

\ms

We introduce the following hypothesis for the involved coefficients in the state equation.

\ms

{\bf(H4)} For $k\ges1$, let
$$A,\bar A_k,\wt A_k,C,\bar C_k,\wt C_k:[0,T]\times\O\to\dbR^{n\times n},\q
B,\bar B_k,\wt B_k,D,\bar D_k,\wt D_k:[0,T]\times\O\to\dbR^{n\times m},$$
be $\dbF$-progressively measurable processes satisfying
%
\bel{Sec4_Condi_ABCD}\left\{\ba{ll}
\ns\ds {\rm(i)}\q\int_0^T\[\esssup_{\o\in\O}|A(s,\o)|+\(\sum_{k\ges1}\dbE|\bar A_k(s)|^2\)^{1\over2}\(\sum_{k\ges1}\dbE|\wt A_k(s)|^2\)^{1\over2}\]ds<\i,\\
\ns\ds
{\rm(ii)}\q\int_0^T\[\esssup_{\o\in\O}|B(s,\o)|^2+\(\sum_{k\ges1}\dbE|\bar
B_k(s)|^2\)
\(\sum_{k\ges1}\dbE|\wt B_k(s)|^2\)\]ds<\i,\\
\ns\ds
{\rm(iii)}\q\int_0^T\[\esssup_{\o\in\O}|C(s,\o)|^2+\(\sum_{k\ges1}\dbE|\bar
C_k(s)|^2\)\(\sum_{k\ges1}
\dbE|\wt C_k(s)|^2\)\]ds<\i,\\
\ns\ds {\rm(iv)}\q \esssup_{s\in [0,T]}\[ \esssup_{\o\in\O}|D(s,\o)|
+\(\sum_{k\ges 1}\dbE |\bar D_k(s)|^2\)^{1\over 2}\(\sum_{k\ges
1}\dbE |\wt D_k(s)|^2\)^{1\over 2}
\]<\i,\ea\right.\ee
and $b(\cd)\in L^2_\dbF(\O;L^1(0,T;\dbR^n))$, $\si(\cd)\in
L^2_\dbF(0,T;\dbR^n)$.

\ms

Note that (H4) ensures the operators $\cA(\cd),\cB(\cd),\cC(\cd),\cD(\cd)$ satisfy (H1).
In fact, take any $\xi\in L^2_{\cF_s}(\O;\dbR^n)$, we have
$$\ba{ll}
\ns\ds\|\cA(s)\xi\|_2=\(\dbE|\cA(s)\xi|^2\)^{1\over2}=\(\dbE\Big|A(s)\xi+\sum_{k\ges1}\bar A_k(s)\dbE[\wt A_k(s)\xi]\Big|^2\)^{1\over2}\\
\ns\ds\les\[\dbE\(|A(s)|^2|\xi|^2\)\]^{1\over2}+\[\dbE\(\sum_{k\ges1}|\bar A_k(s)|\big(\dbE|\wt A_k(s)|^2\big)^{1\over2}\|\xi\|_2\)^2\]^{1\over2}\\
\ns\ds\les\Big\{\esssup_{\o\in\O}|A(s,\o)|+\[\dbE\(\sum_{k\ges1}|\bar A_k(s)|\big(\dbE|\wt A_k(s)|^2\big)^{1\over2}\)^2\]^{1\over2}\Big\}\|\xi\|_2\\
\ns\ds\les\Big\{\esssup_{\o\in\O}|A(s,\o)|+\[\(\dbE\sum_{k\ges1}|\bar
A_k(s)|^2\)\(\sum_{k\ges1} \dbE|\wt
A_k(s)|^2\)\]^{1\over2}\Big\}\|\xi\|_2.\ea$$
Then, \rf{Sec4_Condi_ABCD}-(i) implies that
$$\ba{ll}
\ns\ds\int_0^T\|\cA(s)\|ds\les\int_0^T
\[\esssup_{\o\in\O}|A(s,\o)|+\(\sum_{k\ges1}\dbE|\bar
A_k(s)|^2\)^{1\over2}\(\sum_{k\ges1}\dbE|\wt
A_k(s)|^2\)^{1\over2}\]ds<\i.\ea$$
Therefore, $\cA(\cd)\in \sL^1_\dbF(L^2(\O;\dbR^n))$.

\ms

Similarly, \rf{Sec4_Condi_ABCD}-(ii), (iii), (iv) imply that $\cB(\cd)\in \sL^2_\dbF(L^2(\O;\dbR^m);L^2(\O;\dbR^n))$,
$\cC(\cd)\in \sL^2_\dbF(L^2(\O;\dbR^n))$, $\cD(\cd)\in \sL^\i_\dbF(L^2(\O;\dbR^m);L^2(\O;\dbR^n))$, respectively.

\ms

Next, for the cost functional, we introduce the following hypothesis:

\ms

{\bf(H5)} For $i,j,k\ges1$, let
$$G,\bar G_k,\wt G_k,\h G_{ij}:\O\to\dbR^{n\times n},\q G^\top=G,\q\h G_{ij}^\top=\h G_{ji},$$
be $\cF_T$-measurable and satisfy
%
%
%
\bel{G}\esssup_{\o\in\O}|G(\o)|+\sum_{k\ges1}\dbE\(|\bar G_k|^2+|\wt G_k|^2\)+\sum_{i,j\ges1}|\dbE\h G_{ij}|^2<\i.\ee
and for $i,j,k\ges1$, let
$$\ba{ll}
\ns\ds Q,\bar Q_k,\wt Q_k,\h Q_{ij}:[0,T]\times\O\to\dbR^{n\times n},\q Q(\cd)^\top=Q(\cd),
\q\h Q_{ij}(\cd)^\top=\h Q_{ji}(\cd),\\
\ns\ds R,\bar R_k,\wt R_k,\h R_{ij}:[0,T]\times\O\to\dbR^{m\times m},\q R(\cd)^\top=R(\cd),\q
\h R_{ij}(\cd)^\top=\h R_{ji}(\cd),\\
\ns\ds S,\bar S_k,\h S_{ij}:[0,T]\times\O\to\dbR^{m\times n},\q\wt S_k:[0,T]\times\O\to\dbR^{n\times m},\ea$$
be $\dbF$-progressively measurable processes satisfying
\bel{Sec4_Condi_QRS}\left\{\ba{ll}
\ns\ds{\rm(i)}\q\int_0^T\Big\{\esssup_{\o\in\O}|Q(s,\o)|+\sum_{k\ges1}\dbE\(|\wt Q_k(s)|^2+|\bar Q_k(s)|^2\)\\
\ns\ds\qq\qq+\(\sum_{i,j\ges1}\big|\dbE\h Q_{ij}(s)\big|^2\)^{1\over2}\(\sum_{k\ges1}\dbE|\wt Q_k(s)|^2\)\Big\}ds<\i,\\
\ns\ds{\rm(ii)}\q\esssup_{(s,\o)\in[0,T]\times\O}|R(s,\o)|+\esssup_{s\in[0,T]}\[\sum_{k\ges1}\dbE\(|\wt R_k(s)|^2
+|\bar R_k(s)|^2\)+\sum_{i,j\ges1}\big|\dbE\h R_{ij}(s)\big|^2\]<\i,\\
\ns\ds{\rm(iii)}\q\int_0^T\Big\{\esssup_{\o\in\O}|S(s,\o)|^2+\[\sum_{k\ges1}\dbE\(|\wt Q_k(s)|^2+|\bar S_k(s)|^2
+|\wt S_k(s)|^2+|\wt R_k(s)|^2\)\]^2\\
\ns\ds\qq\qq+\(\sum_{i,j\ges1}\big|\dbE\h
S_{ij}(s)\big|^2\)\(\sum_{k\ges1}\dbE|\wt
Q_k(s)|^2\)\(\sum_{k\ges1}\dbE|\wt
R_k(s)|^2\)\Big\}ds<\i.\ea\right.\ee
Also, $g_0\in L^2(\O;\dbR^n)$, $q_0(\cd)\in
L^2_\dbF(\O;L^1(0,T;\dbR^n)),$ $\rho_0(\cd)\in
L^2_\dbF(0,T;\dbR^m)$, and $\bar q_k:\O\times[0,T]\to\dbR^n$,
$\bar\rho_k:\O\times[0,T]\to\dbR^m$ are $\dbF$-progressively
measurable, and
\bel{Sec4_Condi_q,rho}\sum_{k\ges1}|\dbE\bar g_k|^2+\int_0^T\sum_{k\ges1}\(|\dbE\bar q_k(s)|^2+|\dbE\bar\rho_k(s)|^2\)ds<\i.\ee

\ms
Under (H5), one has all the conditions in (H2) are satisfied. 
Detailed calculations are collected in the appendix.

\ms

Now, for any $(t,x)\in\sD$, our state equation is given by
\bel{Sec4_MFSDE}\left\{\2n\ba{ll}
\ns\ds dX(s) =\Big(A(s)X(s)+\bar\BA(s)^\top \dbE[\wt\BA(s) X(s)]
+B(s)u(s) +\bar \BB(s)^\top\dbE[\wt\BB(s) u(s)]+b(s) \Big) ds\\
\ns\ds \q +\Big( C(s)X(s) +\bar\BC(s)^\top\dbE[\wt\BC(s) X(s)]
+D(s)u(s) +\bar\BD(s)^\top\dbE[\wt\BD(s)u(s)]+\si(s) \Big)
dW(s),\ s\in [t,T],\\
\ns\ds  X(t) = x,\ea\right.\ee
and the quadratic cost functional is
\bel{Sec4_cost}\ba{ll}
\ns\ds J(t,x;u(\cd)) = \dbE\bigg\{ \lan GX(T), X(T) \ran
+2\lan \bar\BG X(T),\dbE[\wt\BG X(T)] \ran +\lan \h\BG
\dbE[\wt\BG X(T)], \dbE[\wt\BG X(T)] \ran\\
\ns\ds\qq +2\lan g_0, X(T) \ran +2\lan \bar\Bg,  \mathbb
E[\wt\BG X(T)] \ran +\int_t^T \Big[ \lan QX,X \ran +2\lan
\bar\BQ X, \dbE[\wt\BQ X] \ran +\lan \h\BQ
\dbE[\wt\BQ X], \dbE[\wt\BQ X] \ran\\
\ns\ds \qq +2\lan SX,u \ran +2\lan \bar\BS X, \dbE[\wt\BR
u] \ran +2\lan \wt\BS \dbE[\wt\BQ X], u \ran +\lan \h\BS \dbE[\wt\BQ X], \dbE[\wt\BR u] \ran
+\lan Ru,u \ran+\lan \dbE[\bar\BR u],\wt\BR u\ran \\
\ns\ds
\qq +\lan \bar\BR u, \dbE[\wt\BR u] \ran +\lan \h\BR \dbE[\wt\BR u], \dbE[\wt\BR u] \ran +2\lan q_0,
X\ran +2\lan \bar\Bq, \dbE[\wt\BQ X] \ran +2\lan
\rho_0, u \ran +2\lan \bar\Brho, \dbE[\wt\BR u] \ran
\Big] ds \bigg\}.\ea\ee
The argument $s$ is suppressed in the above functional.

\ms

In \cite{Buckdahn-Li-Peng 2009}, a well-posedness result of MF-FSDE
was established under Lipschitz condition. In Section 2 of this
paper, we extend the result to a kind of non-Lipschitz case (see
Proposition \ref{Well-posedness-state}) to ensure that, for any
$u(\cd)\in L^2_{\dbF}(t,T;\dbR^m)$, there exists a unique strong
solution $X(\cd)=X(\cd\,;t,x,u(\cd))$ to \rf{Sec4_MFSDE}, and the
cost functional $J(t,x;u(\cd))$ is well-defined. Now, we propose the
MF-LQ stochastic optimal control problem as follows:

\ms

{\bf Problem (MF-LQ).} For any given $(t,x)\in\sD$, find an
admissible control $\bar u(\cd)\in \cU[t,T]$ such that
\bel{} J(t,x;\bar u(\cd)) = \inf_{u(\cd)\in\cU[t,T]}J(t,x;u(\cd)).
\ee
In the above $\bar u(\cd)$ is called an {\it open-loop optimal
control} of Problem (MF-LQ), and the corresponding {\it optimal
state trajectory} $X(\cd;t,x,\bar u(\cdot))$ is denoted by $\bar
X(\cd)$.

\ms

The following result characterizes the optimal control $\bar u(\cd)$
of Problem (MF-LQ).

\bt{Sec4_Th_optimal} \sl Under {\rm(H4)--(H5)}.  For any
$(t,x)\in\sD$ given, $\bar u(\cd)\in\cU[t,T]$ is an open-loop
optimal control of Problem (MF-LQ) at $(t,x)$ with $\bar X(\cd)$
being the corresponding open-loop optimal state process, if and only
if $u(\cd) \mapsto J(t,x;u(\cd))$ is convex and $(\bar X(\cd),\bar
Y(\cd), \bar Z(\cd), \bar u(\cd)) \in
[L^2_\dbF(\O;C([t,T];\dbR^n))]^2 \times L^2_{\dbF}(t,T;\dbR^n)
\times \cU[t,T]$ solves the following system (the argument $s$ is
suppressed):
\bel{Sec4_Hamil_Sys}\left\{\2n\ba{ll}
\ns\ds d\bar X = \Big( A\bar X +\bar\BA^\top \dbE[\wt\BA\bar X]
+B\bar
u +\bar\BB^\top\dbE[\wt\BB\bar u] +b \Big) ds\\
\ns\ds \qq\q +\Big( C\bar X +\bar\BC^\top \dbE[\wt\BC\bar X] +D\bar
u +\bar\BD^\top\dbE[\wt\BD\bar u] +\si \Big) dW,\q s\in
[t,T],\\
\ns\ds d\bar Y = -\Big( A^\top\bar Y +\wt\BA^\top \dbE[\bar\BA\bar
Y] +C^\top\bar Z +\wt\BC^\top \dbE[\bar\BC\bar Z] +Q\bar X +\wt
\BQ^\top \dbE[\bar\BQ\bar X] +\bar\BQ^\top \dbE[\wt\BQ\bar X]\\
\ns\ds \qq\q +\wt\BQ^\top \dbE[\h\BQ]\dbE[\wt\BQ \bar X] +S^\top\bar
u +\bar\BS^\top\dbE[\wt\BR\bar u] +\wt\BQ^\top\dbE[\wt\BS^\top \bar
u] +\wt\BQ^\top \dbE[\h\BS] \dbE[\wt\BR
\bar u] +q \Big) ds\\
\ns\ds\qq\qq\qq\qq+\bar ZdW,\qq s\in [t,T],\\
\ns\ds\bar X(t) =x,\\
\ns\ds\bar Y(T)= G\bar X(T) +\wt\BG^\top\dbE[\bar\BG\bar X(T)]
+\bar\BG^\top\dbE[\wt\BG\bar X(T)] +\wt\BG^\top
\dbE[\h\BG] \dbE[\wt\BG\bar X(T)] +g_0 +\wt\BG^\top \dbE[\bar\Bg],\\
\ns\ds R\bar u +\wt\BR^\top \dbE[\bar\BR\bar u]
+\(\bar\BR^\top+\wt\BR^\top\dbE[\h\BR]\)\dbE[\wt\BR\bar u]
+B^\top\bar Y +\wt\BB^\top \dbE[\bar\BB\bar Y] +D^\top\bar Z +\wt\BD^\top \dbE[\bar\BD\bar Z]\\
\ns\ds \qq +S\bar X +\wt\BR^\top \dbE[\bar\BS\bar X] +\wt \BS\mathbb
E[\wt\BQ\bar X] +\wt\BR^\top \dbE[\h\BS] \dbE[\wt\BQ\bar X] +\rho
=0,\q s\in [t,T].\ea\right.\ee
Moreover, if (H3) holds true, then the above system
\rf{Sec4_Hamil_Sys} admits a unique solution, and Problem (MF-LQ)
admits a unique optimal open-loop control $\bar u(\cd)$. \et

\it Proof. \rm By the above analysis and the calculations collected
in the appendix, we know that (H4)---(H5) imply (H1)---(H2). Then
the results of theorem are obtained by applying Theorems
\ref{optimal-condition} and \ref{Sec5_wellposedness}, and Corollary
\ref{Yu_Sec3.2_Corollary} in Section 3.
\endpf

\ms

\br{} \rm We note that (H3) implies the existence of the inverse
of operator $\cR(\cd)$. In other words, $\bar u(\cd)$ can be solved
from the algebraic equation (the last two lines) in the system
\rf{Sec4_Hamil_Sys}, in terms of $(\bar X(\cd),\bar Y(\cd),\bar Z(\cd))$.
Therefore, system \rf{Sec4_Hamil_Sys} is a coupled mean-field FBSDE.
\er

Assumption (H3) is in the form of operators. In order to suit
mean-field problems, here we try to propose a sufficient condition
of (H3) in the form of matrices. We rewrite the cost functional as
follows:
\bel{Sec4_Sim_Cost}\ba{ll}
\ns\ds J(t,x;u(\cd))=\dbE\bigg\{\lan\BG\wt\BX(T),\wt\BX(T)\ran+2\lan\Bg,\wt\BX(T)\ran
+\int_t^T\[\lan\BQ(s)\BX(s),\BX(s)\ran\\
\ns\ds\qq\qq\qq\qq+2\lan\BS(s)\BX(s),\Bu(s)\ran+\lan\BR(s)\Bu(s),\Bu(s)\ran+2\lan\Bq(s),\BX(s)\ran
+2\lan\Brho(s),\Bu(s)\ran\]ds\bigg\},\ea\ee
where, for any $s\in [t,T]$, the above used notations (introduced in Subsection 2.1) are repeated again as
$$\wt\BX(T)=\begin{pmatrix}X(T)\\ \dbE[\wt\BG X(T)]\end{pmatrix},\qq
\BX(s)=\begin{pmatrix}X(s)\\ \dbE[\wt\BQ(s)X(s)]\end{pmatrix},\qq
\Bu(s)=\begin{pmatrix}u(s)\\ \dbE[\wt\BR(s)u(s)]\end{pmatrix},$$

$$\ba{ll}
\BG = \begin{pmatrix}G&\bar\BG^\top\\ \bar\BG&\h\BG\end{pmatrix}, \q
\BQ(s)=\begin{pmatrix}Q(s) & \bar\BQ(s)^\top \\ \bar\BQ(s) & \h\BQ(s)\end{pmatrix},\q
\BS(s)=\begin{pmatrix}S(s) & \wt\BS(s)\\
\bar\BS(s) & \h\BS(s)\end{pmatrix},
\ea$$

$$
\BR(s) =\begin{pmatrix}R(s) & \bar\BR(s)^\top \\ \bar\BR(s) &
\h\BR(s)\end{pmatrix},\q \Bg = \begin{pmatrix} g_0 \\
\bar\Bg\end{pmatrix},\q \Bq(s)=\begin{pmatrix}q_0(s)\\
\bar\Bq(s)\end{pmatrix}, \q \Brho(s) = \begin{pmatrix}\rho_0(s)\\
\bar\Brho(s)\end{pmatrix}.
$$
We introduce the following condition:
\bel{Sec4_Positive_Assu}
\BG \ges 0, \q
\begin{pmatrix}\BQ(\cd) & \BS(\cd)^\top \\ \BS(\cd) & \BR(\cd)\end{pmatrix}\ges 0, \q
\BR \ges \d \begin{pmatrix}I & 0 \\ 0 & 0 \end{pmatrix}.\ee

\bc{}\sl Condition \eqref{Sec4_Positive_Assu} implies Assumption
(H3). Then under (H4)---(H5) and \eqref{Sec4_Positive_Assu}, all the
results of Theorem {\rm\ref{Sec4_Th_optimal}} hold true.

\ec

\ms

\it Proof. \rm We only need to show \rf{Sec4_Positive_Assu} implies
Assumption (H3). Firstly, for any $\xi\in L^2(\O;\dbR^n)$,
$$\ba{ll}
\ns\ds\dbE\lan \cG\xi,\xi\ran=\dbE\lan
G\xi+\wt\BG^\top\dbE[\bar\BG\xi]+\bar \BG^\top\dbE[\wt\BG\xi]+\wt\BG^\top\dbE[\h \BG]\dbE[\wt\BG\xi], \xi \ran\\
\ns\ds=\dbE\lan
G\xi,\xi\ran+2\dbE\lan\bar\BG\xi,\dbE[\wt\BG\xi]\ran+\dbE\lan\h\BG\dbE[\wt\BG\xi],
\dbE[\wt\BG\xi]\ran=\dbE\lan\BG\wt\Bxi,\wt\Bxi\ran\ges0,\ea$$
where $\wt\Bxi^\top \equiv (\xi^\top, (\dbE[\wt\BG \xi])^\top)$.

\ms

Secondly, for any $s\in [t,T]$, any $\xi\in L^2_{\cF_s}(\O;\dbR^n)$,
any $\eta\in L^2_{\cF_s}(\O;\dbR^m)$, we have (the argument $s$ is
suppressed for simplicity):
$$\ba{ll}
\ns\ds \dbE \left\lan \begin{pmatrix} \cQ &\cS^*\\ \cS &
\cR\end{pmatrix} \begin{pmatrix} \xi\\ \eta \end{pmatrix},
\begin{pmatrix} \xi\\ \eta \end{pmatrix} \right\ran
= \dbE\lan \cQ \xi,\xi \ran+2\dbE\lan \cS\xi,\eta \ran +\dbE\lan\cR\eta,\eta \ran\\
\ns\ds=\Big\{\dbE\lan Q\xi,\xi \ran +2\dbE\lan \bar\BQ\xi, \dbE[\wt\BQ\xi] \ran +\dbE\lan \h\BQ \dbE[\wt\BQ\xi], \dbE[\wt\BQ\xi] \ran\Big\}\\
\ns\ds\q  +2 \Big\{ \dbE\lan S\xi, \eta \ran +\dbE\lan\bar\BS\xi, \dbE[\wt\BR\eta]\ran +\dbE\lan \wt\BS\dbE[\wt\BQ\xi], \eta \ran +\dbE\lan \h\BS\dbE[\wt\BQ\xi], \dbE[\wt\BR\eta] \ran \Big\}\\
\ns\ds\q +\Big\{\dbE\lan R\eta,\eta \ran +2\dbE\lan \bar\BR\eta, \dbE[\wt\BR\eta] \ran +\dbE\lan \h\BR \dbE[\wt\BR\eta], \dbE[\wt\BR\eta] \ran\Big\}\\
\ns\ds= \dbE\lan \BQ\Bxi, \Bxi \ran +2\dbE\lan\BS\Bxi, \BBeta \ran
+\dbE\lan \BR\BBeta, \BBeta \ran=\dbE\left\lan \left(
\begin{array}{ccc} \BQ & \BS^\top\\ \BS & \BR
\end{array} \right) \begin{pmatrix} \Bxi\\ \BBeta \end{pmatrix}, \begin{pmatrix} \Bxi\\ \BBeta \end{pmatrix}
\right\ran \ges 0,\ea$$
where $\Bxi^\top \equiv (\xi^\top, (\dbE[\wt\BQ \xi])^\top )$ and
$\BBeta^\top \equiv (\eta^\top, (\dbE[\wt\BR \eta])^\top)$.

\ms

Thirdly, for any $s\in [t,T]$, any $\eta\in L^2_{\cF_s}(\O;\dbR^m)$,
$$\dbE\lan \cR\eta,\eta \ran -\d\dbE|\eta|^2=\dbE\lan\BR\BBeta,\BBeta\ran-\d
\dbE\left\lan \begin{pmatrix} I & 0\\ 0 & 0 \end{pmatrix}\BBeta,
\BBeta \right\ran
 = \dbE\left\lan \left[ \BR -\d\begin{pmatrix} I & 0\\ 0 & 0 \end{pmatrix} \right] \BBeta, \BBeta \right
 \ran \ges 0.$$
\endpf

\bs

\bex{Sec4_Ex_1} \rm Let the state process $X(\cdot)$ be
one-dimensional, and the terminal cost be
$$
J (t,x;u(\cdot)) = \mbox{Var } (X(T)).
$$
In the form of operators, $J(t,x;u(\cdot)) = \mathbb E \big[
\big(\mathcal G X(T)\big) X(T) \big]$, where
$$
\mathcal G \xi = \xi-\mathbb E[\xi],\qquad \forall\ \xi\in
L^2(\Omega;\mathbb R).
$$
It is easy to check that $\mathcal G\in \mathscr{S}\big(
L^2(\Omega;\dbR) \big)$ and $\mathcal G\ges 0$. However, when
we try to rewrite $J(t,x;u(\cdot))$ in the form of matrices, we find
it admits many representations. For example, one representation is
$$
\begin{aligned}
J(t,x;u(\cdot)) =\ & \mathbb E\Big[ \big( X(T) -\mathbb E[X(T)]
\big)^2 \Big] = \mathbb E\Big[ X(T)^2 -2X(T)\mathbb E[X(T)]
+\big(\mathbb E[X(T)]\big)^2 \Big]\\
=\ & \mathbb E \left[ \left\langle \left( \begin{array}{ccc} 1 & -1\\
-1 & 1
\end{array} \right) \left( \begin{array}{ccc} X(T)\\ \mathbb E[X(T)] \end{array} \right),\  \left( \begin{array}{ccc} X(T)\\ \mathbb E[X(T)] \end{array} \right) \right\rangle \right]
= \mathbb E[\langle \BG_1 \widetilde \BX(T),\ \widetilde \BX(T)
\rangle].
\end{aligned}
$$
We verify that $\BG_1 \ges0$. Another representation is given by
$$
\begin{aligned}
J(t,x;u(\cdot)) =\ & \mathbb E[X(T)^2] - \big(\mathbb E[X(T)]\big)^2
= \mathbb E\Big[ X(T)^2 -\big( \mathbb E[X(T)] \big)^2 \Big]\\
=\ & \mathbb E\left[ \left\langle \left( \begin{array}{ccc} 1 & 0\\
0 & -1\end{array} \right) \left( \begin{array}{ccc} X(T)\\ \mathbb
E[X(T)] \end{array} \right),\ \left( \begin{array}{ccc} X(T)\\
\mathbb E[X(T)] \end{array} \right) \right\rangle \right] = \mathbb
E[\langle \BG_2 \widetilde \BX(T),\ \widetilde \BX(T) \rangle].
\end{aligned}
$$
Obviously, $\BG_2$ is not positive semidefinite. \ex

\br{Sec4.1_Rem_after_Ex}\rm Along the way of the above example, we
can have the following comprehension: a group of operators
satisfying (H3) may be represented by many groups of matrices, and
it is possible that some groups of representation matrices do not
satisfy \eqref{Sec4_Positive_Assu}. This is the reason that we would
like to keep (H3) (instead of \eqref{Sec4_Positive_Assu}) in Theorem
\ref{Sec4_Th_optimal}. Moreover, this comprehension suggests that
there are some advantages to introduce the framework of operators in
the study of mean-field problems. \er

\ms

Next, we try to characterize  the optimal control $\bar u(\cd)$ by the solution of a Fredholm type integral equation of the second kind.

\bp{}\sl Let {\rm(H4)--(H5)} and \rf{Sec4_Positive_Assu} hold. Then for any given $(t,x)\in\sD$, Problem {\rm(MF-LQ)} admits a unique open-loop optimal control $\bar u(\cd)$. Let $(\bar X(\cd),\bar Y(\cd),\bar Z(\cd))$ be the corresponding adapted solution to the FBSDE \rf{Sec4_Hamil_Sys}, and let
$$\left\{
\ba{ll}
\ns\ds\wt\th(s,\o)=-R(s,\o)^{-1}\(B(s,\o)^\top\bar Y(s,\o)+\wt\BB(s,\o)^\top\dbE[\bar\BB(s)\bar Y(s)] +D(s,\o)^\top\bar Z(s,\o)\\
\ns\ds\qq\qq\q+\wt\BD(s,\o)^\top\dbE[\bar\BD(s)\bar Z(s)]+S(s,\o)\bar X(s,\o)+\wt\BR(s,\o)^\top \dbE[\bar\BS(s)\bar X(s)]\\
\ns\ds\qq\qq\q+\wt \BS(s,\o)\dbE[\wt\BQ(s)\bar X(s)]+\wt\BR(s,\o)^\top\dbE[\h\BS(s)]\dbE[\wt\BQ(s)\bar X(s)]+\rho(s,\o)\),\q\o\in\O,~s\in[t,T],\\
\ns\ds\G(s,\o,\o')=-R(s,\o)^{-1}\[\wt\BR(s,\o)^\top\bar\BR(s,\o')+\(\bar\BR(s,\o)^\top+\wt \BR(s,\o)^\top\dbE[\h\BR(s)]\)\wt\BR(s,\o')\],\\
\ns\ds\qq\qq\qq\qq\qq\qq\qq\qq\qq\qq\qq\qq\o,\o'\in\O,~s\in[t,T].\ea\right.$$
Then the following integral equation
\bel{Fredholm}u(s,\o)=\wt\th(s,\o)+\int_\O\G(s,\o,\o')u(s,\o')d\dbP(\o'),\qq\o\in\O,~s\in[t,T],\ee
admits a unique solution which is the open-loop optimal control $\bar u(\cd)$. Moreover,
\bel{bar u}\bar u(s,\o)=\wt\th(s,\o)+\int_\O\F(s,\o,\o')\wt\th(s,\o')d\dbP(\o'),\qq(s,\o)\in[t,T]\times\O,\ee
with $\F(s,\o,\o')$ being the unique solution to the following equation:
\bel{Phi}\F(s,\o,\o')=\G(s,\o,\o')+\int_\O\G(s,\o,\bar\o)\F(s,\bar\o,\o')d\dbP(\bar\o),\qq(s,\o,\o')\in[t,T]\times\O\times\O,\ee

\ep

\it Proof. \rm Note that
$$\dbE\lan\cR u,u\ran=\dbE\bigg\lan\begin{pmatrix}R&\bar \BR^\top\\\bar \BR&\h\BR\end{pmatrix}\begin{pmatrix}
u\\ \dbE[\wt \BR u]\end{pmatrix},\begin{pmatrix}
u\\ \dbE[\wt \BR u]\end{pmatrix}\bigg\ran\ges\d\dbE|u|^2.$$
Thus, $\cR$ is invertible. This means that for any $v\in\cU[t,T]$,
$$\cR u=v$$
admits a unique solution $u\in\cU[t,T]$.

We now write the equality (in \rf{Sec4_Hamil_Sys}) which $\bar u(\cd)$ satisfies  more carefully as follows:
\bel{Sec4_equ_u}
R(t)u(t)+\wt\BR(t)^\top\dbE[\bar\BR(t)u(t)]+\(\bar\BR(t)^\top+\wt\BR(t)^\top\dbE[\h\BR(t)]\)\dbE[\wt\BR(t) u(t)]=\th(t),\ee
with
$$\ba{ll}
\ns\ds\th(t)=-\(B(t)^\top\bar Y(t)+\wt\BB(t)^\top\dbE[\bar\BB(t)\bar Y(t)]+D(t)^\top\bar Z(t) +\wt\BD(t)^\top\dbE[\bar\BD(t)\bar Z(t)]+S(t)\bar X(t)\\
\ns\ds\qq\qq+\wt\BR(t)^\top\dbE[\bar\BS(t)\bar X(t)]+\wt \BS(t)\dbE[\wt\BQ(t)\bar
X(t)]+\wt\BR(t)^\top\dbE[\h\BS(t)]\dbE[\wt\BQ(t)\bar X(t)]+\rho(t)\).\ea$$
Note that \rf{Sec4_equ_u} is indeed a Fredholm type integral equation of the second kind:
$$R(t,\o)u(t,\o)\1n+\3n\int_\O\[\wt\BR(t,\o)^\top\1n\bar\BR(t,\o')\1n+\2n\(\bar\BR(t,\o)^\top\3n+\1n\wt \BR(t,\o)^\top\dbE[\h\BR(t)]\)\wt\BR(t,\o')\]u(t,\o')d\dbP(\o')=\th(t,\o),$$
which can be rewritten into the following standard form of Fredholm integral equation of the second kind:
\bel{Fredholm-standard-Proof}u(t,\o)=\wt\th(t,\o)+\int_\O\G(t,\o,\o')u(t,\o')d\dbP(\o'),\qq\o\in\O,~t\in[0,T].\ee
where
$$\left\{\2n\ba{ll}
\ns\ds\wt\th(t,\o)=R(t,\o)^{-1}\th(t,\o),\\
\ns\ds\G(t,\o,\o')=-R(t,\o)^{-1}\[\wt\BR(t,\o)^\top\bar\BR(t,\o')+\(\bar\BR(t,\o)^\top+\wt \BR(t,\o)^\top\dbE[\h\BR(t)]\)\wt\BR(t,\o')\].\ea\right.$$
Our assumption has guaranteed the well-posedness of the above equation. If we define $\F(t,\o,\o')$ to be the unique solution to the equation \rf{Phi}, and let $\bar u(\cd)$ be defined by \rf{bar u}. Then
$$\ba{ll}
\ns\ds\int_\O\G(s,\o,\o')\bar u(s,\o')d\dbP(\o')=\int_\O\G(s,\o,\o')\[\wt\th(s,\o')+\int_\O\F(s,\o',\bar\o)\wt\th(s,\bar\o)d\dbP(\bar\o)\]
d\dbP(\o')\\
\ns\ds=\int_\O\G(s,\o,\o')\wt\th(s,\o')d\dbP(\o')+\int_\O\[\int_\O\G(s,\o,\o')\F(s,\o',\bar\o)d\dbP(\o')
\]\wt\th(s,\bar\o)d\dbP(\bar\o)\\
\ns\ds=\int_\O\G(s,\o,\o')\wt\th(s,\o')d\dbP(\o')+\int_\O\[\F(s,\o,\o')-\G(s,\o,\o')\]\wt\th(s,\o')
d\dbP(\o')\\
\ns\ds=\int_\O\F(s,\o,\o')\wt\th(s,\o')d\dbP(\o')=\bar u(s,\o)-\wt\th(s,\o).\ea$$
Thus, $\bar u(\cd)$ defined by \rf{bar u} is the solution to the integral equation \rf{Fredholm}. \endpf

\subsection{A generalized mean-variance portfolio selection problem}

In Example \ref{Sec4_Ex_1}, the variance of a random variable is
involved in the cost functional. In this subsection, we shall
investigate a generalized mean-variance portfolio selection problem
arising from mathematical finance, where the cost functional is an
extended version of the one in Example \ref{Sec4_Ex_1}.

\ms

Let us consider a Black-Scholes market model consisting of one
risk-free asset (a bank account or a bond) with interest rate
$r(\cdot)$ and one risky asset (a stock) with appreciation rate
$\mu(\cdot)$ and volatility $\sigma(\cdot)$. We assume that
$r(\cdot)$, $\mu(\cdot)$ and $\sigma(\cdot)$ are nonnegative,
$\mathbb F$-progressively measurable, bounded stochastic processes.
Moreover, there exists a constant $\varepsilon>0$ such that
$\sigma(s)\geq \varepsilon$, for all $s\in [t,T]$. Under the
self-financing assumption, by a standard argument, the wealth
process $X(\cdot)$ of an agent evolves according to the following
FSDE:
\begin{equation}\label{Sec4.1_MV_Sys}
\left\{
\begin{aligned}
& dX(s) = \Big[ r(s)X(s) +\theta(s)\pi(s) \Big] ds
+\pi(s)dW(s),\quad s\in [t,T],\\
& X(t) =x,
\end{aligned}
\right.
\end{equation}
where $x\in L^2_{\mathcal F_t}(\Omega;\mathbb R)$ is the initial
endowment, $\theta(\cdot) = \sigma(\cdot)^{-1}(\mu(\cdot)-r(\cdot))$
is the risk premium, $\pi(\cdot) =
\sigma(\cdot)\widetilde\pi(\cdot)$ and $\widetilde \pi(\cdot)$ is
the dollar amount invested in the stock. $\pi(\cdot)\in \cU[t,T] =
L^2_{\mathbb F}(t,T;\mathbb R)$ is called a portfolio which is the
control process of the agent.

\ms

Let $\beta$ be a given bounded $\mathcal F_T$-measurable random
variable, and there exists a constant $\varepsilon>0$ such that
$\beta\geq\varepsilon$. The cost functional is given by
\begin{equation}\label{Sec4.1_MV_Cost}
J(t,x;\pi(\cdot)) = \mathbb E\Big[ \beta\big( X(T)-\mathbb E[X(T)]
\big)^2 \Big] -2\mathbb E[X(T)],\quad \pi(\cdot)\in \cU[t,T].
\end{equation}
Particularly, when $\beta>0$ is a constant, then the cost functional
is reduced to
$$
J(t,x;\pi(\cdot)) = \beta\mbox{Var }\big(X(T)\big) -2\mathbb
E[X(T)],\quad \pi(\cdot)\in \cU[t,T],
$$
which is a classical mean-variance cost functional studied by many
researchers (see \cite{Andersson-Djehiche 2011, Bjork-Murgoci-Zhou
2014, Pham-Wei 2017} for example). Due to this, for the general
case, we call \eqref{Sec4.1_MV_Cost} a generalized mean-variance
cost functional.

\ms

We pose the following generalized mean-variance portfolio selection
problem.

\ms

{\bf Problem (GMV).} For any given $(t,x)\in\sD$, find a portfolio
$\bar \pi(\cd)\in \cU[t,T]$ such that
\bel{} J(t,x;\bar \pi(\cd)) =
\inf_{\pi(\cd)\in\cU[t,T]}J(t,x;\pi(\cd)). \ee
In the above $\bar \pi(\cd)$ is called an {\it optimal portfolio} of
Problem (GMV), and the corresponding {\it optimal wealth process}
$X(\cd;t,x,\bar \pi(\cdot))$ is denoted by $\bar X(\cd)$.

\ms

Similar to Example \ref{Sec4_Ex_1}, we introduce an operator
$\mathcal G_0 \in \mathscr S\big( L^2(\Omega;\mathbb R) \big)$:
\begin{equation}
\mathcal G_0\xi = \beta\xi -\beta\mathbb E[\xi] -\mathbb E[\beta\xi]
+\mathbb E[\beta]\mathbb E[\xi],\qquad \forall\ \xi\in
L^2(\Omega;\mathbb R),
\end{equation}
as well as an $\mathcal F_T$-measurable random variable $g_0 \equiv
-1$. Then the generalized mean-variance cost functional is rewritten
as $J(t,x;\pi(\cdot)) = \mathbb E\big[\big(\mathcal G_0 X(T)\big)
X(T) +2g_0 X(T)\big]$. For any $\xi \in L^2(\Omega;\mathbb R)$,
$$
\begin{aligned}
& \mathbb E \big[ \big( \mathcal G_0\xi \big) \xi \big] = \mathbb
E\Big\{ \beta\xi^2 -\beta\mathbb E[\xi]\xi -\mathbb E[\beta\xi]\xi
+\mathbb E[\beta]\mathbb E[\xi]\xi \Big\}\\
& = \mathbb E\Big\{ \beta\xi^2 -2\beta\mathbb E[\xi]\xi
+\beta\big(\mathbb E[\xi]\big)^2 \Big\} = \mathbb E\Big\{ \beta
\big( \xi-\mathbb E[\xi] \big)^2 \Big\} \geq 0,
\end{aligned}
$$
i.e., $\mathcal G_0 \geq 0$. However, due to the vanishing of the
control weight, the cost functional defined by
\eqref{Sec4.1_MV_Cost} does not satisfy Assumption (H3). Then the
conclusions in the paper cannot be used directly to solve Problem
(GMV). Next, we try to use the {\it equivalent cost functional
method} introduced in \cite{Yu 2013,Huang-Yu 2014} to overcome the
difficulty.

\ms

It is easy to see that $X(\cdot) -\mathbb E[X(\cdot)]$ satisfies the
following FSDE (suppressing $s$):
\begin{equation}
\left\{
\begin{aligned}
& d\Big( X-\mathbb E[X] \Big) = \Big( rX -\mathbb E[rX] +\theta\pi
-\mathbb E[\theta\pi] \Big) ds +\pi dW,\quad s\in [t,T],\\
& X(t)-\mathbb E[X(t)] = x-\mathbb E[x].
\end{aligned}
\right.
\end{equation}
Let $H(\cdot)$ be a deterministic differential function which will
be determined later. By applying It\^{o}'s formula to $H(\cdot)
(X(\cdot)-\mathbb E[X(\cdot)])^2$, we have
$$
\begin{aligned}
& d \big\{ H(X-\mathbb E[X])^2 \big\}\\
=\ & \Big\{ H'X^2 -2H'\mathbb E[X]X +H'\mathbb E[X]\mathbb E[X]
+2HrX^2 -2H\mathbb E[rX]X -2Hr\mathbb E[X]X +2H\mathbb E[X]\mathbb
E[rX]\\
& +2H\theta X\pi -2HX\mathbb E[\theta\pi] -2H\mathbb E[X]\theta\pi
+2H\mathbb E[X]\mathbb E[\theta\pi] +H\pi^2 \Big\} ds
+2H\big(X-\mathbb E[X] \big)\pi dW.
\end{aligned}
$$
Then,
\begin{equation}
\begin{aligned}
- H(t) \mathbb E\Big[ \big( x-\mathbb E[x] \big)^2 \Big] =\ & - H(T)
\mathbb E\Big[ \big( X(T)-\mathbb E[X(T)] \big)^2
\Big] \\
& + \int_t^T \mathbb E\Big\{ \Big( H'X -H'\mathbb E[X] +2rHX
-H\mathbb E[rX] -Hr\mathbb E[X] \Big) X\\
& +2H\theta \big( X-\mathbb E[X] \big) \pi  +H\pi^2 \Big\} ds.
\end{aligned}
\end{equation}
We notice that, the left hand side of the above is a constant which
is independent of the portfolio $\pi(\cdot)$, and the right hand
side is in a quadratic form of $X(\cdot)$ and $\pi(\cdot)$. Now we
add the constant $-H(t)\mathbb E[(x-\mathbb E[x])^2]$ on the
original cost functional $J(t,x;\cdot)$ to get an equivalent one (in
the sense of evaluating portfolios $\pi(\cdot)$):
\begin{equation}\label{Sec4.1_Equiv_Cost}
\begin{aligned}
J^H(t,x;\pi(\cdot)) =\ & J(t,x;\pi(\cdot)) - H(t) \mathbb E\Big[
\big( x-\mathbb E[x] \big)^2 \Big]\\
=\ & \mathbb E\Big[ \big( \beta-H(T) \big) \big( X(T)-\mathbb
E[X(T)] \big)^2
\Big] -2\mathbb E[X(T)] \\
& + \int_t^T \mathbb E\Big\{ \Big( H'X -H'\mathbb E[X] +2rHX
-H\mathbb E[rX] -Hr\mathbb E[X] \Big) X\\
& +2H\theta \big( X-\mathbb E[X] \big) \pi  +H\pi^2 \Big\} ds.
\end{aligned}
\end{equation}
For the equivalent cost functional $J^H(t,x;\cdot)$, we hope to find
a suitable differential function $H(\cdot)$ such that Assumption
(H3) is satisfied. If in this case, due to the equivalence between
$J(t,x;\cdot)$ and $J^H(t,x;\cdot)$, we can use $J^H(t,x;\cdot)$
instead of $J(t,x;\cdot)$. For this aim, we introduce the following
operators:
\begin{equation}
\left\{
\begin{aligned}
& \mathcal G\xi = \big( \beta-H(T) \big) \xi -\big( \beta-H(T) \big)
\mathbb E[\xi] -\mathbb E \big[ \big( \beta-H(T) \big) \xi \big]
+\mathbb E\big[ \beta-H(T) \big] \mathbb E[\xi],\quad \xi
\in L^2(\Omega;\mathbb R),\\
& \mathcal Q\xi = H'\xi -H'\mathbb E[\xi] +2rH\xi -H\mathbb E[r\xi]
-Hr\mathbb E[\xi],\qquad \xi \in L^2_{\mathcal F_s}(\Omega;\mathbb
R),\quad s\in [t,T],\\
& \mathcal S\xi = H\theta \big( \xi -\mathbb E[\xi] \big),\qquad \xi
\in L^2_{\mathcal F_s}(\Omega;\mathbb
R),\quad s\in [t,T],\\
& \mathcal R\eta = H\eta,\qquad \eta\in L^2_{\mathcal
F_s}(\Omega;\mathbb R),\quad s\in [t,T].
\end{aligned}
\right.
\end{equation}
We also define $g \equiv -1$. Then, \eqref{Sec4.1_Equiv_Cost} can be
rewritten as
\begin{equation}
J^H(t,x;\pi(\cdot)) = \mathbb E\Big[ \mathcal GX(T)\cdot X(T)
+2g\cdot X(T) +\int_t^T \Big( \mathcal QX \cdot X +2\mathcal S X
\cdot \pi +\mathcal R\pi\cdot\pi \Big) ds \Big].
\end{equation}

\ms

Consequently, we have
\begin{equation}
\mathcal S^*\eta = H\big( \theta\eta -\mathbb E[\theta\eta] \big),
\qquad \eta\in L^2_{\mathcal F_s}(\Omega;\mathbb R),\quad s\in
[t,T].
\end{equation}
For any $s\in [t,T]$, any $\xi \in L^2_{\mathcal F_s}(\Omega;\mathbb
R)$, we calculate
\begin{equation}
\begin{aligned}
& \mathbb E\Big[ \big( \mathcal Q-\mathcal S^*\mathcal
R^{-1}\mathcal S \big) \xi\cdot \xi \Big]\\
=\ & H' \mathbb E\Big[ \big(\xi-\mathbb E[\xi]\big)\xi \Big]
+H\mathbb E\Big[(2r-\theta^2)\big(\xi-\mathbb E[\xi]\big) \xi\Big]
-H\mathbb E\Big[ \theta^2\big(\xi-\mathbb E[\xi]\big) \Big] \mathbb
E[\xi]\\
=\ & \mathbb E \Big[ \big( H' +H(2r-\theta^2) \big) \big(
\xi-\mathbb E[\xi] \big)^2 \Big] +H \mathbb E\Big[ (2r-2\theta^2)
\big( \xi-\mathbb E[\xi] \big) \Big] \mathbb E[\xi].
\end{aligned}
\end{equation}
We introduce the following assumption.

\ms

{\bf(H6).} The process $r(\cd)-\th(\cd)^2$ is deterministic.

\medskip

\noindent Under Assumption (H6),
\begin{equation}
\mathbb E\Big[ \big( \mathcal Q-\mathcal S^*\mathcal R^{-1}\mathcal
S \big) \xi\cdot \xi \Big] = \mathbb E \Big[ \big( H'
+H(2r-\theta^2) \big) \big( \xi-\mathbb E[\xi] \big)^2 \Big].
\end{equation}
Now, we give a sufficient condition of (H3) as follows.
\begin{equation}\label{Sec4.1_H3}
\left\{
\begin{aligned}
\mbox{(i).} \quad & H(T) \leq \beta,\\
\mbox{(ii).} \quad & H(s) >0,\quad s\in [t,T],\\
\mbox{(iii).} \quad & H'(s) +H(s)\big(2r(s)-\theta(s)^2\big) \geq
0,\quad s\in [t,T].
\end{aligned}
\right.
\end{equation}
Since $r(\cdot)$ and $\theta(\cdot)$ are bounded, then there exists
a constant $K>0$ such that $|2r(s)-\theta(s)^2| \leq K$ for all
$s\in [t,T]$. We define
\begin{equation}
H(s) = \varepsilon e^{-K(T-s)},\quad s\in [t,T].
\end{equation}
It is easy to verify that $H(\cdot)$ defined above satisfies all the
requirements in \eqref{Sec4.1_H3}.

\bp{}\sl Under Assumption (H6), Problem (GMV) admits a unique
optimal portfolio $\bar\pi(\cdot)\in \mathcal U[t,T]$, where
$(X(\cdot),\bar\pi(\cdot),\bar Y(\cdot))\in L^2_{\mathbb
F}(\Omega;C([t,T];\mathbb R)) \times \mathcal U[t,T] \times
L^2_{\mathbb F}(\Omega;C([t,T];\mathbb R))$ is the unique solution
to the following FBSDE:
\begin{equation}\label{Sec4.1_MV_FBSDE}
\left\{
\begin{aligned}
& d\bar X(s) = \Big[ r(s)\bar X(s) +\theta(s)\bar\pi(s) \Big] ds
+\bar \pi(s) dW(s),\quad s\in [t,T],\\
& d\bar Y(s) = -r(s)\bar Y(s)ds -\theta(s)\bar Y(s)dW(s),\quad s\in [t,T],\\
& \bar X(t) =x,\quad \bar Y(T) = \beta \bar X(T) -\beta\mathbb
E[\bar X(T)] -\mathbb E[\beta X(T)] +\mathbb E [\beta ] \mathbb
E[\bar X(T)] -1.
\end{aligned}
\right.
\end{equation}

\ep

We notice that FBSDE \eqref{Sec4.1_MV_FBSDE} is not in the classical
form.

\ms

\it Proof. \rm Firstly, under Assumption (H6), for the optimal
control problem \eqref{Sec4.1_MV_Sys} and \eqref{Sec4.1_Equiv_Cost},
Assumptions (H3)---(H5) are satisfied. Then Theorem
\ref{Sec4_Th_optimal} works to ensure that, the optimal control
problem \eqref{Sec4.1_MV_Sys} and \eqref{Sec4.1_Equiv_Cost} admits a
unique optimal control $\tilde \pi(\cdot) \in \mathcal U[t,T]$,
where $(\tilde X(\cdot),\tilde \pi(\cdot),\tilde Y(\cdot), \tilde
Z(\cdot)) \in L^2_{\mathbb F}(\Omega;C([t,T];\mathbb R)) \times
\mathcal U[t,T] \times L^2_{\mathbb F}(\Omega;C([t,T];\mathbb R))
\times L^2_{\mathbb F}(t,T;\mathbb R)$ is the unique solution to the
following system:
\begin{equation}\label{Sec4.1_Equivalent_Ham_Sys}
\left\{
\begin{aligned}
& H\tilde \pi +\theta\tilde Y +\tilde Z +H\theta\Big( \tilde
X-\mathbb E [\tilde X]\Big) =0, \quad
s\in [t,T],\\
& d\tilde X = \Big[ r\tilde X +\theta\tilde\pi \Big] ds
+\tilde \pi dW,\quad s\in [t,T],\\
& d\tilde Y = -\Big\{r\tilde Y +H'\tilde X-H'\mathbb E[\tilde X]
+2rH\tilde X -H\mathbb E[r\tilde X] -Hr\mathbb E[\tilde X]\\
& \qquad\qquad +H\theta\tilde\pi
-H\mathbb E[\theta\tilde\pi] \Big\} ds +\tilde Z dW,\quad s\in [t,T],\\
& \tilde X(t) =x,\\
& \tilde Y(T) = \big( \beta -H(T) \big)\tilde X(T)
-\big(\beta-H(T)\big)\mathbb E[\tilde X(T)] -\mathbb E\big[
\big(\beta-H(T)\big)\tilde X(T) \big]\\
& \qquad\qquad +\mathbb E\big[ \beta -H(T) \big] \mathbb E[\tilde
X(T)] -1.
\end{aligned}
\right.
\end{equation}
Since the cost functionals $J^H(t,x;\cdot)$ and $J(t,x;\cdot)$ are
equivalent, then $\bar \pi(\cdot) = \tilde \pi(\cdot)$ is also the
unique optimal portfolio for Problem (GMV). Consequently, $\bar
X(\cdot) = \tilde X(\cdot)$ is the corresponding optimal wealth
process.

\ms

Secondly, by Theorem \ref{optimal-condition}, the optimal pair
$(\bar X(\cdot), \bar \pi(\cdot))$ of Problem (GMV) as well as the
corresponding pair of adjoint processes $(\bar Y(\cdot),\bar
Z(\cdot)) \in L^2_{\mathbb F}(\Omega;C([t,T];\mathbb R)) \times
L^2_{\mathbb F}(t,T;\mathbb R)$ satisfies the following system:
\begin{equation}\label{Sec4.1_MV_Ham_Sys}
\left\{
\begin{aligned}
& \theta(s) \bar Y(s) +\bar Z(s) =0,\quad s\in [t,T],\\
& d\bar X(s) = \Big[ r(s)\bar X(s) +\theta(s)\bar\pi(s) \Big] ds
+\bar \pi(s) dW(s),\quad s\in [t,T],\\
& d\bar Y(s) = -r(s)\bar Y(s)ds +\bar Z(s)dW(s),\quad s\in [t,T],\\
& \bar X(t) =x,\quad \bar Y(T) = \beta \bar X(T) -\beta\mathbb
E[\bar X(T)] -\mathbb E[\beta X(T)] +\mathbb E [\beta ] \mathbb
E[\bar X(T)] -1.
\end{aligned}
\right.
\end{equation}
Moreover, we can verify that there exists a relationship between
systems \eqref{Sec4.1_Equivalent_Ham_Sys} and
\eqref{Sec4.1_MV_Ham_Sys}:
\begin{equation}\label{Sec4.1_Transform_Ham_Sys}
\tilde X = \bar X,\quad \tilde \pi =\bar \pi, \quad \tilde Y = \bar
Y -H\big( \bar X-\mathbb E[\bar X] \big),\quad \tilde Z =\bar Z
-H\bar \pi,\quad s\in [t,T].
\end{equation}
We notice that the above transformation is invertible, which implies
the uniqueness of \eqref{Sec4.1_MV_Ham_Sys}. Now we have proved that
system \eqref{Sec4.1_MV_Ham_Sys} admits a unique solution.

\ms

Finally, by substituting $\bar Z(\cdot) = -\theta(\cdot)\bar
Y(\cdot)$ into the forward and backward SDEs of
\eqref{Sec4.1_MV_Ham_Sys}, we complete the proof.
\endpf

\br{} \rm As the same as the classical stochastic LQ problems with
matrix coefficients, (H3) is regarded as a {\it regular} (or
standard) assumption. But the (generalized) mean-variance portfolio
selection problem is a typical example in the non-regular case. In
the analysis of this subsection, we introduced Assumption (H6). We
hope we can get rid of it in our future research.

\er

\section{Conclusion}

It is the mean-field LQ Problem (MF-LQ) that inspires us to study
the LQ optimal control problem with operator coefficients (i.e.,
Problem (OLQ)). As we known, \rf{Sec4_MFSDE}--\rf{Sec4_cost} is a new form of mean-field LQ problems.
Besides, all the coefficients are allowed to be random in our study.
As a start, we only study the open-loop case.  The closed-loop case
of the control problems, as well as differential
games are under our investigation.

\section*{Acknowledgement}

This work was carried out during the stay of Qingmeng Wei and
Zhiyong Yu at University of Central Florida, USA. They would like to
thank the hospitality of Department of Mathematics, and the
financial support from China Scholarship Council. Also, the authors would like to thank the anonymous referees for their careful suggestive comments.

\section{Appendix.}

We collect some detailed computations of Section 4 below. Recall that
\bel{}\left\{\2n\ba{ll}
\ns\ds\cG\xi=G\xi+\wt\BG^\top\dbE[\bar\BG\xi]+\bar\BG^\top\dbE[\wt\BG\xi]
+\wt\BG^\top\dbE[\h\BG]\dbE[\wt\BG\xi],\q\xi\in L^2(\O;\dbR^n),\\
\ns\ds\cQ(s)\xi=Q(s)\xi+\wt\BQ(s)^\top\dbE[\bar\BQ(s)\xi]+\bar\BQ(s)^\top\dbE[\wt\BQ(s)\xi]
+\wt\BQ(s)^\top\dbE[\h\BQ(s)]\dbE[\wt\BQ(s)\xi],\q\xi\in
L^2_{\cF_s}(\O;\dbR^n),\\
\ns\ds\cS(s)\xi=S(s)\xi+\wt\BR(s)^\top\dbE[\bar\BS(s)\xi]+\wt\BS(s)^\top\dbE[\wt\BQ(s)\xi]
+\wt\BR(s)^\top\dbE[\h\BS(s)]\dbE[\wt\BQ(s)\xi],\q\xi\in
L^2_{\cF_s}(\O;\dbR^n),\\
\ns\ds\cR(s)\eta=R(s)\eta+\wt\BR(s)^\top\dbE[\bar\BR(s)\eta]+\bar\BR(s)^\top\dbE[\wt\BR(s)\eta]
+\wt\BR(s)^\top\dbE[\h\BR(s)]\dbE[\wt\BR(s)\eta],\q\eta\in
L^2_{\cF_s}(\O;\dbR^m),\\
\ns\ds g=g_0+\wt\BG^\top\dbE[\bar\Bg],\q
q(s)=q_0(s)+\wt\BQ(s)^\top\dbE[\bar\Bq(s)],\q
\rho(s)=\rho_0(s)+\wt\BR(s)^\top\dbE[\bar\Brho(s)].\ea\right.\ee
Then
$$\ba{ll}
\ns\ds\|\cG\xi\|_2\les\|G\xi\|_2+\|\wt\BG^\top\dbE[\bar\BG\xi]\|_2+\|\bar\BG^\top\dbE[\wt\BG\xi]\|_2
+\|\wt\BG^\top\dbE[\h\BG]\dbE[\wt\BG\xi]\|_2\\
\ns\ds\les\(\dbE|G\xi|^2\)^{1\over2}+\(\dbE\[\sum_{k\ges1}|\bar
G_k|\big(\dbE|\wt
G_k|^2\big)^{1\over2}\|\xi\|_2\]^2\)^{1\over2}+\(\dbE\[\sum_{k\ges1}|\wt
G_k|
\big(\dbE|\bar G_k|^2\big)^{1\over2}\|\xi\|_2\]^2\)^{1\over2}\\
\ns\ds\qq\qq\q +\(\dbE\[\sum_{i,j\ges1}|\wt G_i| |\dbE\h G_{ij} |\big(\dbE|\wt G_j|^2\big)^{1\over2}\|\xi\|_2\]^2\)^{1\over2}\\
\ns\ds\les\Big\{\esssup_{\o\in\O}|G(\o)|+\(\sum_{k\ges1}\dbE|\bar
G_k|^2\)^{1\over2}\(\sum_{k\ges1}
\dbE|\wt G_k|^2\)^{1\over2}+\(\sum_{k\ges1}\dbE|\wt G_k|^2\)^{1\over2}\(\sum_{k\ges1}\dbE|\bar G_k|^2\)^{1\over2}\\
\ns\ds\qq\qq\q+\(\sum_{i\ges1}\dbE|\wt G_i|^2 \)^{1\over2} \(\sum_{i\ges1}\[\sum_{j\ges1} |\dbE\h G_{ij}|\big(\dbE|\wt G_j|^2\big)^{1\over2} \]^2\)^{1\over2}\Big\}\|\xi\|_2\\
\ns\ds\les\Big\{\esssup_{\o\in\O}|G(\o)|+\sum_{k\ges1}\dbE|\bar
G_k|^2+\sum_{k\ges1}\dbE|\wt G_k|^2+\(\sum_{i,j\ges1}|\dbE\h
G_{ij}|^2\)^{1\over2}\(\sum_{k\ges1}\dbE|\wt
G_k|^2\)\Big\}\|\xi\|_2.\ea$$
This implies that
$$\|\cG\|\les\esssup_{\o\in\O}|G(\o)|+\sum_{k\ges1}\dbE|\bar G_k|^2+\sum_{k\ges1}\dbE|\wt G_k|^2+\(\sum_{i,j\ges1}|\dbE\h G_{ij}|^2\)^{1\over2}\(\sum_{k\ges1}\dbE|\wt G_k|^2\).$$
With the same idea, we have
$$\ba{ll}
\ns\ds\|\cQ(s)\xi\|_2\les\|Q(s)\xi\|_2+\|\wt\BQ(s)^\top\dbE[\bar\BQ(s)\xi]\|_2+\|\bar\BQ(s)^\top
\dbE[\wt\BQ(s)\xi]\|_2+\|\wt\BQ(s)^\top\dbE[\h\BQ(s)]\dbE[\wt\BQ(s)\xi]\|_2\\
\ns\ds=\(\dbE|Q(s)\xi|^2\)^{1\over2}+\(\dbE\big|\sum_{k\ges1}\wt Q_k(s)^\top\dbE[\bar Q_k(s)\xi]\big|^2\)^{1\over2}+\(\dbE\big|\sum_{k\ges1}\bar Q_k(s)^\top\dbE[\wt Q_k(s)\xi]\big|^2\)^{1\over2}\\
\ns\ds\qq+\(\dbE\big|\sum_{i,j\ges1}\wt Q_i(s)^\top\dbE[\h Q_{ij}(s)]\dbE[\wt Q_j(s)\xi]\big|^2\)^{1\over2}\\
\ns\ds\les\esssup_{\o\in\O}|Q(s,\o)|\|\xi\|_2+\[\dbE\(\sum_{k\ges1}|\wt Q_k(s)|\big(\dbE|\bar Q_k(s)|^2\big)^{1\over2}\|\xi\|_2\)^2\]^{1\over2}+\[\dbE\(\sum_{k\ges1}|\bar Q_k(s)|\big(\dbE|\wt Q_k(s)|^2\big)^{1\over2}\|\xi\|_2\)^2\]^{1\over2}\\
\ns\ds\qq+\Big\{\dbE\[\(\sum_{i\ges1}|\wt Q_i(s)|^2\)\sum_{i\ges1}\(\sum_{j\ges1}|\dbE\h Q_{ij}(s)|\big|\dbE[\wt Q_j(s)\xi]\big|\)^2\]\Big\}^{1\over2}\\
\ns\ds\les\Big\{\esssup_{\o\in\O}|Q(s,\o)|+2\(\sum_{k\ges1}\dbE|\wt Q_k(s)|^2\)^{1\over2}
\(\sum_{k\ges1}\dbE|\bar Q_k(s)|^2\)^{1\over2}\\
\ns\ds\qq+\(\sum_{i\ges1}\dbE|\wt Q_i(s)|^2\)^{1\over2}\[\sum_{i\ges1}\(\sum_{j\ges1}|\dbE\h Q_{ij}(s)|\big(\dbE|\wt Q_j(s)|^2\big)^{1\over2}\)^2\]^{1\over2}\Big\}\|\xi\|_2\\
\ns\ds\les\Big\{\esssup_{\o\in\O}|Q(s,\o)|+\sum_{k\ges1}\dbE|\wt
Q_k(s)|^2+\sum_{k\ges1}\dbE|\bar Q_k(s)|^2+\(\sum_{i,j\ges1}|\dbE\h
Q_{ij}(s)|^2\)^{1\over2}\(\sum_{k\ges1}\dbE|\wt
Q_k(s)|^2\)\Big\}\|\xi\|_2,\ea$$
which leads to
$$\|\cQ(s)\|\les\esssup_{\o\in\O}|Q(s,\o)|+\sum_{k\ges1}\dbE|\wt Q_k(s)|^2+\sum_{k\ges1}\dbE|\bar Q_k(s)|^2+\(\sum_{i,j\ges1}|\dbE\h
Q_{ij}(s)|^2\)^{1\over2}\(\sum_{k\ges1}\dbE|\wt Q_k(s)|^2\).$$
Next,
$$\ba{ll}
\ns\ds\|\cR(s)\eta\|_2\les\|R(s)\eta\|_2+\|\wt\BR(s)^\top\dbE[\bar\BR(s)\eta]\|_2+\|\bar\BR(s)^\top
\dbE[\wt\BR(s)\eta]\|_2+\|\wt\BR(s)^\top\dbE[\h\BR(s)]\dbE[\wt\BR(s)\eta]\|_2\\
\ns\ds=\(\dbE|R(s)\xi|^2\)^{1\over2}+\(\dbE\big|\sum_{k\ges1}\wt R_k(s)^\top\dbE[\bar R_k(s)\eta]\big|^2\)^{1\over2}+\(\dbE\big|\sum_{k\ges1}\bar R_k(s)^\top\dbE[\wt R_k(s)\eta]\big|^2\)^{1\over2}\\
\ns\ds\qq+\(\dbE\big|\sum_{i,j\ges1}\wt R_i(s)^\top\dbE[\h R_{ij}(s)]\dbE[\wt R_j(s)\eta]\big|^2\)^{1\over2}\\
\ns\ds\les\esssup_{\o\in\O}|R(s,\o)|\|\eta\|_2+\[\dbE\(\sum_{k\ges1}|\wt R_k(s)|\big(\dbE|\bar R_k(s)|^2\big)^{1\over2}\|\eta\|_2\)^2\]^{1\over2}+\[\dbE\(\sum_{k\ges1}|\bar R_k(s)|\big(\dbE|\wt R_k(s)|^2\big)^{1\over2}\|\eta\|_2\)^2\]^{1\over2}\\
\ns\ds\qq+\Big\{\dbE\[\(\sum_{i\ges1}|\wt R_i(s)|^2\)\sum_{i\ges1}\(\sum_{j\ges1}|\dbE\h R_{ij}(s)| \big|\dbE[\wt R_j(s)\eta]\big|\)^2\]\Big\}^{1\over2}\\
\ns\ds\les\Big\{\esssup_{\o\in\O}|R(s,\o)|+2\(\sum_{k\ges1}\dbE|\wt R_k(s)|^2\)^{1\over2}
\(\sum_{k\ges1}\dbE|\bar R_k(s)|^2\)^{1\over2}\\
\ns\ds\qq+\(\sum_{i\ges1}\dbE|\wt R_i(s)|^2\)^{1\over2}\[\sum_{i\ges1}\(\sum_{j\ges1}|\dbE\h R_{ij}(s)| \big(\dbE|\wt R_j(s)|^2\big)^{1\over2}\)^2\]^{1\over2}\Big\}\|\eta\|_2\\
\ns\ds\les\Big\{\esssup_{\o\in\O}|R(s,\o)|+\sum_{k\ges1}\dbE|\wt
R_k(s)|^2+\sum_{k\ges1}\dbE|\bar R_k(s)|^2+\(\sum_{i,j\ges1}|\dbE\h
R_{ij}(s)|^2\)^{1\over2}\(\sum_{k\ges1}\dbE|\wt
R_k(s)|^2\)\Big\}\|\eta\|_2,\ea$$
which leads to
$$\|\cR(s)\|\les\esssup_{\o\in\O}|R(s,\o)|+\sum_{k\ges1}\dbE|\wt R_k(s)|^2+\sum_{k\ges1}\dbE|\bar R_k(s)|^2+\(\sum_{i,j\ges1}|\dbE\h R_{ij}(s)|^2\)^{1\over2}\(\sum_{k\ges1}\dbE|\wt R_k(s)|^2\).$$
Further,
$$\ba{ll}
\ns\ds\|\cS(s)\xi\|_2\les\|S(s)\xi\|_2+\|\wt\BR(s)^\top\dbE[\bar\BS(s)\xi]\|_2+\|\wt\BS(s)^\top
\dbE[\wt\BQ(s)\xi]\|_2+\|\wt\BR(s)^\top\dbE[\h\BS(s)]\dbE[\wt\BQ(s)\xi]\|_2\\
\ns\ds=\(\dbE|S(s)\xi|^2\)^{1\over2}+\(\dbE\big|\sum_{k\ges1}\wt R_k(s)^\top\dbE[\bar S_k(s)\xi]\big|^2\)^{1\over2}+\(\dbE\big|\sum_{k\ges1}\wt S_k(s)^\top\dbE[\wt Q_k(s)\xi]\big|^2\)^{1\over2}\\
\ns\ds\qq+\(\dbE\big|\sum_{i,j\ges1}\wt R_i(s)^\top\dbE[\h S_{ij}(s)]\dbE[\wt Q_j(s)\xi]\big|^2\)^{1\over2}\\
\ns\ds\les\esssup_{\o\in\O}|S(s,\o)|\|\xi\|_2+\[\dbE\(\sum_{k\ges1}|\wt R_k(s)|\big(\dbE|\bar S_k(s)|^2\big)^{1\over2}\|\xi\|_2\)^2\]^{1\over2}+\[\dbE\(\sum_{k\ges1}|\wt S_k(s)|\big(\dbE|\wt Q_k(s)|^2\big)^{1\over2}\|\xi\|_2\)^2\]^{1\over2}\\
\ns\ds\qq+\Big\{\dbE\[\(\sum_{i\ges1}|\wt R_i(s)|^2\)\sum_{i\ges1}\(\sum_{j\ges1} |\dbE\h S_{ij}(s)| \big|\dbE[\wt Q_j(s)\xi]\big|\)^2\]\Big\}^{1\over2}\\
\ns\ds\les\Big\{\esssup_{\o\in\O}|S(s,\o)|+\(\sum_{k\ges1}\dbE|\wt R_k(s)|^2\)^{1\over2}
\(\sum_{k\ges1}\dbE|\bar S_k(s)|^2\)^{1\over2}+\(\sum_{k\ges1}\dbE|\wt S_k(s)|^2\)^{1\over2}
\(\sum_{k\ges1}\dbE|\wt Q_k(s)|^2\)^{1\over2}\\
\ns\ds\qq+\(\sum_{i\ges1}\dbE|\wt R_i(s)|^2\)^{1\over2}\[\sum_{i\ges1}\(\sum_{j\ges1} |\dbE\h S_{ij}(s)| \big(\dbE|\wt Q_j(s)|^2\big)^{1\over2}\)^2\]^{1\over2}\Big\}\|\xi\|_2\\
\ns\ds\les\Big\{\esssup_{\o\in\O}|S(s,\o)|+\sum_{k\ges1}\dbE|\wt R_k(s)|^2+\sum_{k\ges1}\dbE|\bar S_k(s)|^2+\sum_{k\ges1}\dbE|\wt S_k(s)|^2+
\sum_{k\ges1}\dbE|\wt Q_k(s)|^2\\
\ns\ds\qq+\(\sum_{i,j\ges1}\big|\dbE\h
S_{ij}(s)\big|^2\)^{1\over2}\(\sum_{k\ges1}\dbE|\wt
R_k(s)|^2\)^{1\over2}\(\sum_{k\ges1}\dbE|\wt
Q_k(s)|^2\)^{1\over2}\Big\}\|\xi\|_2,\ea$$
which leads to
$$\ba{ll}
\ns\ds\|\cS(s)\|\les\esssup_{\o\in\O}|S(s,\o)|+\sum_{k\ges1}\dbE|\wt R_k(s)|^2+\sum_{k\ges1}\dbE|\bar S_k(s)|^2+\sum_{k\ges1}\dbE|\wt S_k(s)|^2+
\sum_{k\ges1}\dbE|\wt Q_k(s)|^2\\
\ns\ds\qq\qq\qq+\(\sum_{i,j\ges1}\big|\dbE\h
S_{ij}(s)\big|^2\)^{1\over2}\(\sum_{k\ges1}\dbE|\wt
R_k(s)|^2\)^{1\over2}\(\sum_{k\ges1}\dbE|\wt
Q_k(s)|^2\)^{1\over2}.\ea$$
Finally,
$$\ba{ll}
\ns\ds\|g\|_2\les\|g_0\|_2+\|\wt\BG^\top\dbE[\bar\Bg]\|_2=\|g_0\|_2+\[\dbE\Big|\sum_{k\ges1}\wt G_k^\top\dbE[\bar g_k]\Big|^2\]^{1\over2}\\
\ns\ds\qq\les\|g_0\|_2+\(\sum_{k\ges1}\dbE|\wt
G_k|^2\)^{1\over2}\(\sum_{k\ges1}|\dbE\bar g_k|^2\)^{1\over2}
\les\|g_0\|_2+\sum_{k\ges1}\dbE|\wt G_k|^2+\sum_{k\ges1}\big|\dbE\bar g_k\big|^2;\\
\ns\ds \Big\|\int_0^T |q(s)|ds\Big\|_2 \les \Big\| \int_0^T
|q_0(s)|ds \Big\|_2 +\int_0^T \| \wt\BQ(s)^\top\dbE[\bar\Bq(s)]
\|_2ds\\
\ns\ds\qq\q = \Big\| \int_0^T |q_0(s)|ds \Big\|_2 +\int_0^T \[\dbE\Big|\sum_{k\ges1}\wt Q_k(s)^\top\dbE[\bar q_k(s)]\Big|^2\]^{1\over2} ds\\
\ns\ds\qq\q\les\Big\| \int_0^T |q_0(s)|ds \Big\|_2+ \int_0^T
\(\sum_{k\ges1}\dbE|\wt
Q_k(s)|^2\)^{1\over2}\(\sum_{k\ges1}|\dbE\bar q_k(s)|^2\)^{1\over2}
ds\\
\ns\ds\qq\q
\les\Big\| \int_0^T |q_0(s)|ds \Big\|_2+ \int_0^T\(\sum_{k\ges1}\dbE|\wt Q_k(s)|^2\)ds +\int_0^T \(\sum_{k\ges1}|\dbE \bar q_k(s)|^2\) ds;\\
\ns\ds\|\rho(s)\|_2\les\|\rho_0(s)\|_2+\|\wt\BR(s)^\top\dbE[\bar\Brho(s)]\|_2
=\|\rho_0(s)\|_2+\[\dbE\Big|\sum_{k\ges1}\wt R_k(s)^\top\dbE[\bar\rho_k(s)]\Big|^2\]^{1\over2}\\
\ns\ds\qq\q\les\|\rho_0(s)\|_2+\(\sum_{k\ges1}\dbE|\wt
R_k(s)|^2\)^{1\over2}\(\sum_{k\ges1}|\dbE\bar\rho_k(s)|^2\)^{1\over2}
\les\|\rho_0(s)\|_2+\sum_{k\ges1}\dbE|\wt
R_k(s)|^2+\sum_{k\ges1}\big|\dbE\bar\rho_k(s)\big|^2.\ea$$

\end{document}